\documentclass[10pt]{article}
\usepackage{lmodern}
\usepackage{amsmath}
\usepackage[T1]{fontenc}
\usepackage[utf8]{inputenc}
\usepackage{authblk}
\usepackage{amsfonts}
\usepackage{graphicx}
\usepackage{rotating}
\usepackage{amssymb}
\usepackage[english]{babel}
\usepackage{color}
\usepackage{amsthm}
\usepackage{graphicx}
\usepackage{mathrsfs}
\usepackage{makecell}
\usepackage{microtype}
\usepackage{enumerate}
\usepackage{mathscinet}
\usepackage{array}
\usepackage{multirow}
\usepackage{booktabs}
\usepackage{enumerate}
\usepackage[cal=boondoxo,bb=ams]{mathalfa}
\usepackage{hyperref}
\hypersetup{hidelinks}
\newtheorem{defn}{Definition}[section]

\newtheorem{theorem}{Theorem}[section]

\newtheorem{coro}{Corollary}[section]
\newtheorem{remark}{Remark}[section]
\newtheorem{exam}{Example}[section]

\newcommand{\ml}{\mathcal}
\newcommand{\mb}{\mathbb}

\DeclareMathOperator{\non}{non}
\DeclareMathOperator{\lin}{lin}

\title{Interplay effects on blow-up of weakly coupled systems for semilinear wave equations with general nonlinear memory terms}
\author[1]{Wenhui Chen\thanks{Corresponding author: Wenhui Chen (wenhui.chen.math@gmail.com)}}
\affil[1]{Institute of Applied Analysis, Faculty of Mathematics and Computer Science, Technical University Bergakademie Freiberg, Pr\"uferstra{\ss}e 9, 09596 Freiberg, Germany}

\date{}
\begin{document}

\maketitle
\begin{abstract}
In this paper, we study weakly coupled systems for semilinear wave equations with distinct nonlinear memory terms in general forms, and the corresponding single semilinear equation with general nonlinear memory terms. Thanks to Banach's fixed point theorem, we prove local (in time) existence of solutions with the $L^1$ assumption on the memory kernels. Then, blow-up results for energy solutions are derived applying iteration methods associated with slicing procedure. We investigate interactions on the blow-up conditions under different decreasing assumptions on the memory term. Particularly, a new threshold for the kernels on interplay effect is found. Additionally, we give some applications of our results on semilinear wave equations and acoustic wave equations.\\
	
	
	\noindent\textbf{AMS Classification (2010)} Primary: 35L71, 35B44; Secondary: 35L05, 35B33.
\end{abstract}
\fontsize{12}{15}
\selectfont

\section{Introduction}\setcounter{equation}{0}
In this paper, we mainly work on the Cauchy problem of weakly coupled systems for semilinear wave equations with distinct nonlinear memory terms in general ways, namely,
\begin{align}\label{Eq Semi Wave Distinct Memory}
\begin{cases}
u_{tt}-\Delta u=g_1\ast|v|^p,&x\in\mb{R}^n,\ t>0,\\
v_{tt}-\Delta v=g_2\ast|u|^q,&x\in\mb{R}^n,\ t>0,\\
(u,u_t,v,v_t)(0,x)=(u_0,u_1,v_0,v_1)(x),&x\in\mb{R}^n,
\end{cases}
\end{align}
where $p,q>1$, and the convolution nonlinear terms with respect to time variable are denoted by
\begin{align*}
(g_k\ast|w|^r)(t,x):=\int_0^tg_k(t-\tau)|w(\tau,x)|^r\mathrm{d}\tau,
\end{align*}
with the time-dependent memory kernels or the so-called relaxation functions $g_k=g_k(t)>0$ for $k=1,2$. In the above notation of the convolution, we denote $w=u,v$ and $r=p,q$. Our main attempt of the paper is to investigate blow-up conditions for the weakly coupled systems \eqref{Eq Semi Wave Distinct Memory} influenced by various assumptions on the memory kernels. Moreover, as the special case of \eqref{Eq Semi Wave Distinct Memory}, blow-up results for the corresponding single semilinear wave equation with memory nonlinearities to \eqref{Eq Semi Wave Distinct Memory} also will be shown.

Let us present a historical overview on some results for semilinear wave equations, which are strongly linked with our model and the motivation of considering the weakly coupled systems \eqref{Eq Semi Wave Distinct Memory}. In the last forty years, the following Cauchy problem for semilinear wave equations:
\begin{align}\label{Eq Single semilinear wave}
\begin{cases}
u_{tt}-\Delta u=f(u;p),&x\in\mb{R}^n,\ t>0,\\
(u,u_t)(0,x)=(u_0,u_1)(x),&x\in\mb{R}^n,
\end{cases}
\end{align}
with $p>1$ has caught a lot of attention. Concerning the semilinear wave equation with the power nonlinearity, i.e \eqref{Eq Single semilinear wave} with $f(u;p)=|u|^p$, the critical exponent is the so-called Strauss exponent $p_{\mathrm{Str}}(n)$, which is the positive root of the quadratic equation $(n-1)p^2_{\mathrm{Str}}-(n+1)p_{\mathrm{Str}}-2=0$ for $n\geqslant 2$, to be specific, it can be represented explicitly by
\begin{align*}
p_{\mathrm{Str}}(n):=\frac{n+1+\sqrt{n^2+10n-7}}{2(n-1)}\ \ \mbox{for}\ \ n\geqslant2,
\end{align*}
and in the one spatial dimensional case $p_{\mathrm{Str}}(1):=\infty$. Note that the critical exponent describes the threshold between global (in time) existence of small data weak solutions and blow-up of local (in time) small data weak solutions. For the finding of the Strauss' conjecture and its proof on the critical exponent, we refer to the classical works \cite{Joh79,Kato80,Str81,Sid84,Sch85,Gla81-g,Gla81-b,Zhou95,LS96,GLS97,Tataru2001,Jiao03,YordanovZhang2006,Zhou07}. The previous results, expressly, proved blow-up of weak solutions provided $1<p\leqslant p_{\mathrm{Str}}(n)$. Let us now turn to the study of wave equations with nonlinear memory terms. Recently, the authors of \cite{ChenPalmieri20} considered the Riemann-Liouville fractional integral of $1-\gamma$ of the $p$ power of solution as the nonlinear terms on the semilinear wave equation \eqref{Eq Single semilinear wave}, exactly,
\begin{align}\label{Riemann-Liouville nonlinearity}
f(u;p)=\frac{1}{\Gamma(1-\gamma)}\int_0^t(t-\tau)^{-\gamma}|u(\tau,x)|^p\mathrm{d}\tau\ \ \mbox{with}\ \ \gamma\in(0,1),
\end{align}
where $\Gamma$ stands for the Euler integral of the second kind. They investigated a generalized Strauss exponent $p_0(n,\gamma)$ defined by
\begin{align*}
p_0(n,\gamma):=\frac{n+3-2\gamma+\sqrt{n^2+(14-4\gamma)n+(3-2\gamma)^2-8}}{2(n-1)}\ \ \mbox{for} \ \ n\geqslant2,
\end{align*}
and in the one spatial dimensional case $p_0(1,\gamma):=\infty$. Actually, $p_0(n,\gamma)$ is the positive root of the quadratic equation $(n-1)p^2_0-(n+3-2\gamma)p_0-2=0$. The authors of \cite{ChenPalmieri20} proved blow-up of energy solutions to the semilinear wave equation \eqref{Eq Single semilinear wave} with nonlinearities \eqref{Riemann-Liouville nonlinearity} if $1<p\leqslant p_0(n,\gamma)$ for $n\geqslant2$, and $p>1$ for $n=1$, which satisfy $\lim_{\gamma\to 1^-}p_0(n,\gamma)=p_{\mathrm{Str}}(n)$ for all $n\geqslant 1$. For this reason, it seems reasonable to consider the general nonlinear memory terms $f(u;p)=g\ast|u|^p$ by taking suitable assumptions on the memory kernel $g=g(t)$. However, not only global (in time) existence but also blow-up results are still unknown. It seems interesting to find the influence of the memory kernel on these results. We will answer this question in Corollary \ref{Coro subcritical} by considering the special case of \eqref{Eq Semi Wave Distinct Memory} from the point of view of blow-up for energy solutions. Concerning other studies of hyperbolic equations with nonlinear memory terms \eqref{Riemann-Liouville nonlinearity}, we refer the interested readers to \cite{Fino2011,D'Abbicco2014NoDEA,DAbb14,Lai-Liu-Zhao-2017,Dannawi-Kirane-Fino,DaoFino2020}.

On the other hand, it is known that the critical curve in the $p-q$ plane for the Cauchy problem for weakly coupled systems of semilinear wave equations
\begin{align}\label{Eq Semi Wave Coupled}
\begin{cases}
u_{tt}-\Delta u=|v|^p,&x\in\mb{R}^n,\ t>0,\\
v_{tt}-\Delta v=|u|^q,&x\in\mb{R}^n,\ t>0,\\
(u,u_t,v,v_t)(0,x)=(u_0,u_1,v_0,v_1)(x),&x\in\mb{R}^n,
\end{cases}
\end{align}
is given by the cubic relation
\begin{align*}
\alpha_{\mathrm{W}}(p,q):=\max\left\{\frac{p+2+q^{-1}}{pq-1},\frac{q+2+p^{-1}}{pq-1}\right\}=\frac{n-1}{2}.
\end{align*}
Concerning the work on this critical curve, we refer to \cite{DelSanto-Georgiev-Mitidieri1997,DelSanto1997,DelSanto-Mitidieri1998,Agemi-Kurokawa-Takamura2000,Kurokawa-Takamura2003,Kurokawa2005,Georgiev-Takamura-Zhou2006,Said-Kirane2011,Kurokawa-Takamura-Wakasa2012}. Particularly, every local (in time) weak solution to the weakly coupled system \eqref{Eq Semi Wave Coupled} blows up if $\alpha_{\mathrm{W}}(p,q)\geqslant(n-1)/2$. We should remark that the research of semilinear weakly coupled systems is not just a trivial generalization of those for single semilinear equations, specifically in the case when $p\neq q$. More detail expansions on this effect were pointed out in the recent paper \cite{PalTak19}.

The main goal of this paper is to investigate blow-up results for the weakly hyperbolic coupled systems \eqref{Eq Semi Wave Distinct Memory} with distinct and general memory kernels. Specifically, we are interested in the interplay effects between $g_1(t)$ and $g_2(t)$ on the blow-up conditions.

 Roughly speaking, as we will show in Section \ref{Section Main result}, when the memory kernel decreases slower than $t^{-1}$ (slow decrease), we may feel the influence from the memory kernel on the condition for blow-up of energy solutions. Nonetheless, this phenomenon will disappear when the memory kernel decreases faster than $t^{-1}$ (fast decrease), and in this case the blow-up condition for the exponents $p,q$ is the same as those for the weakly coupled systems \eqref{Eq Semi Wave Coupled}. Therefore, one of our novelties is to derive a new threshold $t^{-1}$ for the long time behavior of memory kernels in the weakly coupled systems \eqref{Eq Semi Wave Distinct Memory} by using iteration methods. Additionally, in the case when both memory kernels decrease slower than $t^{-1}$, up to the author best knowledge, the blow-up condition with combined effects from $g_1(t)$ and $g_2(t)$ is attractive. We will give some applications of the obtained results which are the wave equations with the Riemann-Liouville fractional integral type nonlinear memory terms as the example for slow decay memory kernels, and the semilinear Moore-Gibson-Thompson equations in the conservative case as the instance for fast decay memory kernels.
 
 Lastly, we should mention that the study of blow-up for \eqref{Eq Semi Wave Distinct Memory} is not a simple generalization of those in the previous researches on semilinear wave equations with power nonlinearities. For one thing, due to the memory terms on the nonlinearities, it seems that the classical coupled Kato's type lemma in \cite{DelSanto-Georgiev-Mitidieri1997} does not work well in our model. For another, we now do not take any explicit forms of nonlinear memory terms, which give us some technical difficulties in the treatments. To overcome these difficulties, in the present paper, we will employ iteration methods, especially, when the memory kernel decreases fast, slicing procedure (see, for example, \cite{Agemi-Kurokawa-Takamura2000,ChenPalmieri201901}) will be used. 
 
 The remaining part of the paper is organized as follows. We will show our main results on local (in time) existence and blow-up, moreover, give some applications on these results in Section \ref{Section Main result}. We next prove existence of local (in time) solution by applying Banach's fixed point theorem and $L^2-L^2$ estimates in Section \ref{Section Proof of LOCAL EXISTENCE}. Then, the proof of the blow-up results by using iteration methods will be shown in Section \ref{Section Proof of THM 1}. Eventually, final remarks concerning some open problems in Section \ref{Section final remark} complete the present paper.
 
 \medskip
 
 \noindent\textbf{Notation: } We give some notations to be used in this paper. We write $f\lesssim g$ when there exists a positive constant $C$ such that $f\leqslant Cg$ and, similarly, for $f\gtrsim g$. The relation $f\asymp g$ means that $f\lesssim g\lesssim f$. We denote $\lceil r\rceil:=\min\{C\in\mb{Z}: r\leqslant C\}$ as the ceiling function. Moreover, $B_R$ denotes the ball around the origin with radius $R$ in $\mathbb{R}^n$.
 
\section{Main results}\label{Section Main result}\setcounter{equation}{0}
Let us first state the result for local (in time) existence of solutions of the weakly coupled systems \eqref{Eq Semi Wave Distinct Memory} under a assumption on the memory kernels.
\begin{theorem}\label{Thm Local existence}
Let $(u_0,u_1,v_0,v_1)\in(H^1(\mb{R}^n)\times L^2(\mb{R}^n))\times (H^1(\mb{R}^n)\times L^2(\mb{R}^n))$ compactly support in a ball with some radials $R>0$.  We assume $p,q>1$ such that $p,q\leqslant n/(n-2)$ when $n\geqslant3$. Then, there exists a positive $T$ and a uniquely determined local (in time) mild solution
\begin{align*}
(u,v)\in\left(\ml{C}([0,T],H^1(\mb{R}^n))\cap \ml{C}^1([0,T],L^2(\mb{R}^n))\right)^2
\end{align*}
to the weakly coupled systems \eqref{Eq Semi Wave Distinct Memory} with $g_1(t),g_2(t)\in L^1([0,T])$ satisfying $\mathrm{supp}\,u(t,\cdot),\,\mathrm{supp}\,v(t,\cdot)\subset B_{R+t}$ for any $t\in[0,T]$.
\end{theorem}

Before stating our main results for blow-up of solutions, motivated by \cite{Lai-Takamura-Wakasa2017} we first define a suitable energy solutions of the weakly coupled systems \eqref{Eq Semi Wave Distinct Memory}.
\begin{defn}\label{Defn Energy Solution}
Let $(u_0,u_1,v_0,v_1)\in(H^1(\mb{R}^n)\times L^2(\mb{R}^n))\times (H^1(\mb{R}^n)\times L^2(\mb{R}^n))$. We say that $(u,v)$ is an energy solution of the hyperbolic coupled systems \eqref{Eq Semi Wave Distinct Memory} on $[0,T)$ if
\begin{align*}
&u\in\ml{C}([0,T),H^1(\mb{R}^n))\cap \ml{C}^1([0,T),L^2(\mb{R}^n))\ \ \mbox{with}\ \ g_2\ast|u|^q\in L_{\mathrm{loc}}^1([0,T)\times \mb{R}^n),\\
&v\in\ml{C}([0,T),H^1(\mb{R}^n))\cap \ml{C}^1([0,T),L^2(\mb{R}^n))\ \ \mbox{with}\ \ g_1\ast|v|^p\in L_{\mathrm{loc}}^1([0,T)\times \mb{R}^n),
\end{align*}
fulfills $(u,v)(0,\cdot)=(u_0,v_0)$ in $H^1(\mb{R}^n)\times H^1(\mb{R}^n)$ and the integral relations
\begin{align}\label{Eq Defn Energy 1}
&\int_0^t\int_{\mb{R}^n}(-u_t(s,x)\phi_s(s,x)+\nabla u(s,x)\cdot\nabla\phi(s,x))\mathrm{d}x\mathrm{d}s+\int_{\mb{R}^n}u_t(t,x)\phi(t,x)\mathrm{d}x\notag\\
&=\int_{\mb{R}^n}u_1(x)\phi(0,x)\mathrm{d}x+\int_0^t\int_{\mb{R}^n}(g_1\ast|v|^p)(s,x)\phi(s,x)\mathrm{d}x\mathrm{d}s
\end{align}
as well as 
\begin{align}\label{Eq Defn Energy 2}
&\int_0^t\int_{\mb{R}^n}(-v_t(s,x)\psi_s(s,x)+\nabla v(s,x)\cdot\nabla\psi(s,x))\mathrm{d}x\mathrm{d}s+\int_{\mb{R}^n}v_t(t,x)\psi(t,x)\mathrm{d}x\notag\\
&=\int_{\mb{R}^n}v_1(x)\psi(0,x)\mathrm{d}x+\int_0^t\int_{\mb{R}^n}(g_2\ast|u|^q)(s,x)\psi(s,x)\mathrm{d}x\mathrm{d}s
\end{align}
for any test functions $\phi,\psi\in\ml{C}_0^{\infty}([0,T)\times\mb{R}^n)$ and any $t\in(0,T)$.
\end{defn}

By applying once integration by parts in \eqref{Eq Defn Energy 1} and \eqref{Eq Defn Energy 2}, one may immediately observe the resultant equalities
\begin{align*}
&\int_0^t\int_{\mb{R}^n}u(s,x)(\phi_{ss}(s,x)-\Delta\phi(s,x))\mathrm{d}x\mathrm{d}s+\int_{\mb{R}^n}(u_t(t,x)\phi(t,x)-u(t,x)\phi_s(t,x))\mathrm{d}x\\
&=\int_{\mb{R}^n}(u_1(x)\phi(0,x)-u_0(x)\phi_s(0,x))\mathrm{d}x+\int_0^t\int_{\mb{R}^n}(g_1\ast|v|^p)(s,x)\phi(s,x)\mathrm{d}x\mathrm{d}s
\end{align*}
as well as
\begin{align*}
&\int_0^t\int_{\mb{R}^n}v(s,x)(\psi_{ss}(s,x)-\Delta\psi(s,x))\mathrm{d}x\mathrm{d}s+\int_{\mb{R}^n}(v_t(t,x)\psi(t,x)-v(t,x)\psi_s(t,x))\mathrm{d}x\\
&=\int_{\mb{R}^n}(v_1(x)\psi(0,x)-v_0(x)\psi_s(0,x))\mathrm{d}x+\int_0^t\int_{\mb{R}^n}(g_2\ast|u|^q)(s,x)\psi(s,x)\mathrm{d}x\mathrm{d}s
\end{align*}
 satisfying the definition of weak solutions of the hyperbolic system \eqref{Eq Semi Wave Distinct Memory} by letting $t\to T$.


To describe the interplay effects of $g_1(t)$ and $g_2(t)$, we now introduce a monotonously increasing function by
\begin{align*}
\ml{L}(t):=\underbrace{\ln(\ln(\cdots\ln(}_{r+1\,\,\text{times}\,\,\ln}\underbrace{\exp(\exp(\cdots\exp}_{r\,\,\text{times}\,\,\exp}(1)))+t)))
\end{align*}
for any $r\in\mb{N}_0$ and $t\geqslant0$. It carries the initial value for $t=0$ such that $\ml{L}(0)=0$. The aim of constructing $\ml{L}(t)$ is to find a nonnegative function goes to $\infty$ as $t\to\infty$ with the slower increasing rate. Later, we may choose arbitrary nonnegative integer $r$ in all results, even sufficiently large $r$.

\begin{theorem}\label{Thm subcritical}
Let $p,q>1$ such that $p,q\leqslant n/(n-2)$ when $n\geqslant 3$. Let us consider $g_k(t)\in L^1([0,T])$ for $k=1,2$ such that one of the following condition holds:
\begin{itemize}
	\item if $g_k(t)\gtrsim t^{-1}$ for any $t\geqslant t_0$ with $t_0\in[0,T)$ and $g'_k(t)\leqslant0$ for $k=1,2$, then we consider $(p,q)$ so that the next estimate is satisfied:
\begin{align}\label{Hypothesis Condition 1}
g_1(t)g_2(t)\max\left\{(g_1(t))^{q-1}t^{2q+1/p},(g_2(t))^{p-1}t^{2p+1/q}\right\}\gtrsim t^{\frac{n-1}{2}(pq-1)-3}\ml{L}(t);
\end{align}
\item if $g_k(t)\lesssim t^{-1}$ for any $t\geqslant t_0$ with $t_0\in[0,T)$ and $g_k(t)\in\ml{C}^2([0,T])$ with $g_k''(0)>0$ for $k=1,2$, then we consider $(p,q)$ so that the next estimate is satisfied:
\begin{align}\label{Hypothesis Condition 2}
\max\left\{\frac{p+2+q^{-1}}{pq-1},\frac{q+2+p^{-1}}{pq-1}\right\}>\frac{n-1}{2}.
\end{align}
\end{itemize}
 Let us suppose that $(u_0,u_1,v_0,v_1)\in(H^1(\mb{R}^n)\times L^2(\mb{R}^n))\times(H^1(\mb{R}^n)\times L^2(\mb{R}^n))$ are nonnegative and compactly supported functions with supports contained in $B_R$ for some $R>0$ such that $u_1,v_1$ are not identically zero. Let $(u,v)$ be the local (in time) energy solution to the weakly coupled system \eqref{Eq Semi Wave Distinct Memory} according to Definition \ref{Defn Energy Solution}. Then, these solutions fulfill
\begin{align*}
\mathrm{supp}\,u,\,\mathrm{supp}\,v\subset\{(t,x)\in[0,T)\times\mb{R}^n:|x|\leqslant R+t\},
\end{align*}
and the solution $(u,v)$ blows up in finite time, i.e. $T<\infty$.
\end{theorem}
\begin{remark}\label{Rem Sup 01}
	Indeed, the assumption  $g_k'(t)\leqslant0$ can be replaced by: $g_k(t)\gtrsim \tilde{g}_k(t)>0$ for all $t\geqslant0$, where $\tilde{g}'_k(t)\leqslant 0$. One may see the explanation in Remark \ref{Remark General} later. Thus, it can solve a more general case even it has oscillations. For instance, we can choose $g_k(t)=(3+2\sin t)t^{-\gamma}$ with $\gamma\in[0,1)$ so that $g_k(t)\geqslant \tilde{g}_k(t)=t^{-\gamma}$. Here, we restricted $\gamma<1$ since the local existence result holds if $g_k(t)\in L^1([0,T])$.
\end{remark}
\begin{remark}\label{Rem sup 02}
	We may observe from Theorem \ref{Thm subcritical} that the decay function $t^{-1}$ somehow is the threshold of memory kernels on the blow-up result. Precisely, providing that the memory kernels decay slower than $t^{-1}$, we may notice the interplay effect between $g_1(t)$ as well as $g_2(t)$ in the blow-up condition \eqref{Hypothesis Condition 1}. Nevertheless, if the memory kernels decay faster than $t^{-1}$, then they do not influence on the blow-up condition \eqref{Hypothesis Condition 2} any more, which coincides with those for weakly coupled systems \eqref{Eq Semi Wave Coupled}. We should underline that the condition \eqref{Hypothesis Condition 1} continuous to \eqref{Hypothesis Condition 2} as $g_k(t)\to t^{-1}$ for $k=1,2$. What's more, we have obtained the result not only for $g_k(t)\gtrsim t^{-1}$ but also for $g_k(t)\lesssim t^{-1}$, which means that the blow-up result stated in Theorem \ref{Thm subcritical} covers a general class of memory kernels belonging to $ L^1([0,T])$.
\end{remark}
\begin{remark}
	 Let $p,q>1$ such that $p,q\leqslant n/(n-2)$ when $n\geqslant 3$. Concerning the case that only one memory kernel decreases slower than $t^{-1}$, namely,
	\begin{align*}
	\min\{g_1(t),g_2(t)\}\lesssim t^{-1}\lesssim \max\{g_1(t),g_2(t)\}\ \ \mbox{for any}\ \ t\geqslant t_0,
	\end{align*}
	where $t_0\in[0,T)$,  we believe the energy solutions $(u,v)$ of the weakly coupled systems \eqref{Eq Semi Wave Distinct Memory} blows up in finite time, i.e. $T<\infty$, if $g_k(t)\in L^1([0,T])$ for $k=1,2$ such that
\begin{itemize}
	\item  $g_1(t)\in\ml{C}^2([0,T])$ with $g''_1(0)>0$, $g_2'(t)\leqslant0$  and we consider $(p,q)$ so that the next estimate is satisfied:
	\begin{align*}
	g_2(t)\max\left\{t^{1-q+2q+1/p},(g_2(t))^{p-1}t^{2p+1/q}\right\}\gtrsim t^{\frac{n-1}{2}(pq-1)-2}\ml{L}(t),
	\end{align*}
	when $g_1(t)\lesssim t^{-1}\lesssim g_2(t)$ for any $t\geqslant t_0$;
	\item $g_2(t)\in\ml{C}^2([0,T])$ with $g''_2(0)>0$, $g_1'(t)\leqslant0$ and we consider $(p,q)$ so that the next estimate is satisfied:
	\begin{align*}
	g_1(t)\max\left\{(g_1(t))^{q-1}t^{2q+1/p},t^{1-p+2p+1/q}\right\}\gtrsim t^{\frac{n-1}{2}(pq-1)-2}\ml{L}(t),
	\end{align*}
	when $g_2(t)\lesssim t^{-1}\lesssim g_1(t)$ for any $t\geqslant t_0$;
\end{itemize}
	under the assumption that $(u_0,u_1,v_0,v_1)\in(H^1(\mb{R}^n)\times L^2(\mb{R}^n))\times(H^1(\mb{R}^n)\times L^2(\mb{R}^n))$ are nonnegative and compactly supported functions with supports contained in $B_R$ for some $R>0$ such that $u_1,v_1$ are not identically zero. We expect this conjecture can be proved by combining the proof of Case 1 and Case 2 in Section \ref{Section Proof of THM 1} without any additional technical difficulties. Precisely, to treat the nonlinear term containing $\min\{g_1(t),g_2(t)\}$, one may directly repeat the step of those in Case 1. On the other hand, to deal with the nonlinear term containing $\max\{g_1(t),g_2(t)\}$, one may follow the analogous procedure to those stated in Case 2.
\end{remark}
\begin{remark}\label{Rem Sup 03}
The construction of the function $\ml{L}(t)$ on the right-hand side of the condition \eqref{Hypothesis Condition 1} is used to guarantee
\begin{align}\label{Remark condition 1}
g_1(t)g_2(t)\max\left\{(g_1(t))^{q-1}t^{2q+1/p},(g_2(t))^{p-1}t^{2p+1/q}\right\} t^{-\frac{n-1}{2}(pq-1)+3}\geqslant C_{\mathrm{Suit}}
\end{align}
for any $t\geqslant t_{\mathrm{Suit}}$, where $t_{\mathrm{Suit}}\geqslant0$ is a suitable constant and $C_{\mathrm{Suit}}=C_{\mathrm{Suit}}(u_0,u_1,v_0,v_1,n,p,R,p,q)$ is a suitable positive constant. Indeed, due to the fact that by fixing the size of initial data, dimension, radius of the support of initial data and the power exponents of the nonlinearities, i.e. we determined the constant $C_{\mathrm{Suit}}$, there exists a nonnegative constant $t_{\mathrm{Suit}}$ such that
\begin{align*}
\ml{L}(t)\geqslant C_{\mathrm{Suit}}(n,p,R,p,q)
\end{align*}
for any $t\geqslant t_{\mathrm{Suit}}$, where we used the property that $\lim_{t\to\infty}\ml{L}(t)=\infty$.
\end{remark}
\begin{remark}
Indeed, the condition $g_k(t)\in L^1([0,T])$ for the case when $g_k(t)\lesssim t^{-1}$ for $t\geqslant t_0$  with $t_0\in[0,T)$ is trivial. Because of the continuity of $g_k(t)$ for $t\in[0,\max\{t_0,\tilde{t}_0\}]$ with $\tilde{t}_0\in(0,T)$, we obtain $g_k(t)\lesssim 1$ for any $t\in[0,\max\{t_0,\tilde{t}_0\}]$. Thus, 
\begin{align*}
\int_0^Tg_k(t)\mathrm{d}t=\int_0^{\max\{t_0,\tilde{t}_0\}}g_k(t)\mathrm{d}t+\int_{\max\{t_0,\tilde{t}_0\}}^Tg_k(t)\mathrm{d}t\lesssim \max\{t_0,\tilde{t}_0\}+\ln(T)-\ln (\max\{t_0,\tilde{t}_0\}).
\end{align*}
In other words, one derives $g_k(t)\in L^1([0,T])$ for any $k=1,2$.
\end{remark}

As we known in the previous studies on the weakly coupled systems and single semilinear equations, it is clear that in the special case when the weakly coupled systems carrying $p=q$, the result is wholly symmetric to what occurs in the case of a single semilinear equation. Then, we can get the next result by taking $p=q$ and $g(t)=g_1(t)=g_2(t)$ in the last theorem, which also can be proved by strictly following the procedure of those for Theorem \ref{Thm subcritical}. What's more, similarly to Theorem \ref{Thm Local existence}, under the assumption $g(t)\in L^1([0,T])$, the local (in time) solution to the single semilinear equation uniquely exists if $p>1$ when $n=1,2$ and $1<p\leqslant n/(n-2)$ when $n\geqslant 3$.

\begin{coro}\label{Coro subcritical}
Let $p>1$ such that $p\leqslant n/(n-2)$ when $n\geqslant 3$. Let us consider $g(t)\in L^1([0,T])$ such that one of the following condition holds:
	\begin{itemize}
		\item if $g(t)\gtrsim t^{-1}$ for any $t\geqslant t_0$ with $t_0\in[0,T)$ and $g'(t)\leqslant0$, then we consider $(p,q)$ so that the next estimate is satisfied:
		\begin{align}\label{Coro Hypothesis Condition 1}
		g(t)\gtrsim t^{\frac{n-1}{2}(p-1)-2-\frac{1}{p}}\ml{L}(t);
		\end{align}
		\item if $g(t)\lesssim t^{-1}$ for any $t\geqslant t_0$ with $t_0\in[0,T)$ and $g(t)\in\ml{C}^2([0,T])$ with $g''(0)>0$, then we consider $(p,q)$ so that the next estimate is satisfied:
		\begin{align}\label{Coro Hypothesis Condition 2}
		1<p<p_{\mathrm{Str}}(n),
		\end{align}
		where $p_{\mathrm{Str}}(n)$ denotes the Strauss exponent and it has been defined in the introduction.
	\end{itemize}
	Let us suppose that $(u_0,u_1)\in H^1(\mb{R}^n)\times L^2(\mb{R}^n)$ are nonnegative and compactly supported functions with supports contained in $B_R$ for some $R>0$ such that $u_1$ is not identically zero. Let $u$ be the local (in time) energy solution to 
\begin{align}\label{Eq Single semilinear wave Memory}
\begin{cases}
u_{tt}-\Delta u=g\ast|u|^p,&x\in\mb{R}^n,\ t>0,\\
(u,u_t)(0,x)=(u_0,u_1)(x),&x\in\mb{R}^n.
\end{cases}
\end{align}
	Then, these solutions fulfill
	\begin{align*}
	\mathrm{supp}\,u\subset\{(t,x)\in[0,T)\times\mb{R}^n:|x|\leqslant R+t\},
	\end{align*}
	and the solution $u$ blows up in finite time, i.e. $T<\infty$.
\end{coro}
\begin{remark}
Actually, the local (in time) energy solution of \eqref{Eq Single semilinear wave Memory} has the similar definition to those in Definition \ref{Defn Energy Solution} such that
\begin{align*}
u\in\ml{C}([0,T),H^1(\mb{R}^n))\cap \ml{C}^1([0,T),L^2(\mb{R}^n))\ \ \mbox{with}\ \ g\ast|u|^p\in L_{\mathrm{loc}}^1([0,T)\times \mb{R}^n)
\end{align*}
fulfills $u(0,\cdot)=u_0$ in $H^1(\mb{R}^n)$ and the integral relation
\begin{align*}
&\int_0^t\int_{\mb{R}^n}(-u_t(s,x)\phi_s(s,x)+\nabla u(s,x)\cdot\nabla\phi(s,x))\mathrm{d}x\mathrm{d}s+\int_{\mb{R}^n}u_t(t,x)\phi(t,x)\mathrm{d}x\notag\\
&=\int_{\mb{R}^n}u_1(x)\phi(0,x)\mathrm{d}x+\int_0^t\int_{\mb{R}^n}(g\ast|u|^p)(s,x)\phi(s,x)\mathrm{d}x\mathrm{d}s
\end{align*}
for any test function $\phi\in\ml{C}_0^{\infty}([0,T)\times\mb{R}^n)$ and any $t\in(0,T)$.
\end{remark}
\begin{remark}
We should point out that if $g_1(t),g_2(t)$ are replaced by $g(t)$, then  Remarks \ref{Rem Sup 01}, \ref{Rem sup 02} and \ref{Rem Sup 03} also hold for Corollary \ref{Coro subcritical}.
\end{remark}

\subsection{Application of Corollary \ref{Coro subcritical}}
In this subsection, we will give some applications (examples) of Corollary \ref{Coro subcritical} by taking some special kinds of memory kernels. These applications can be coincided with the developed results in \cite{ChenPalmieri201901,ChenPalmieri20}. Throughout this part, we assume $p>1$ such that $p\leqslant n/(n-2)$ when $n\geqslant 3$.
\begin{exam}
	Let us consider the semilinear wave equation \eqref{Eq Single semilinear wave Memory} with the Riemann-Liouville fractional integrals of order $1-\gamma$ and $\gamma\in(0,1)$, that is the nonlinear term \eqref{Riemann-Liouville nonlinearity}. Obviously, it holds that $g(t)\gtrsim t^{-1}$ for any $t\geqslant 0$ and $g(t)\in L^1([0,T])$. We now employ Corollary \ref{Coro subcritical} with the condition \eqref{Coro Hypothesis Condition 1} leading to
	\begin{align*}
	(n-1)p^2-(n+3-2\gamma)p-2<0.
	\end{align*}
	In other words, let us assume initial data satisfying the same assumptions as those in Corollary \ref{Coro subcritical} and $1<p<p_0(n,\gamma)$, where $p_0(n,\gamma)$ stands for a generalized Strauss exponent defined in the introduction. Then, every energy solution of the semilinear wave equation \eqref{Eq Single semilinear wave} carrying \eqref{Riemann-Liouville nonlinearity} blows up in finite time. Thus, this result coincides with Theorem 1 in \cite{ChenPalmieri20}.
\end{exam}
\begin{exam}\label{Examp single 2}
Let us consider the semilinear wave equation \eqref{Eq Single semilinear wave Memory} with exponential decay memory kernel $g(t)=\mathrm{e}^{-t/\beta}$ with $\beta>0$. Note that $g(t)\lesssim t^{-1}$ for $t\geqslant t_0(\beta)$. From Corollary \ref{Coro subcritical}, the blow-up condition \eqref{Coro Hypothesis Condition 2}, exactly, is $1<p<p_{\mathrm{Str}}(n)$. Thanks to the derivative relation for the exponential decay function that
\begin{align*}
\partial_t(g\ast|u|^p)(t,x)=|u(t,x)|^p-\frac{1}{\beta}(g\ast|u|^p)(t,x),
\end{align*}
 we may transfer the Cauchy problem in this example into a special kind of semilinear Moore-Gibson-Thompson equations in the conservative case describing acoustic waves (one may the detail introduction in the recent paper \cite{Pellicer-Said2019}), namely,
\begin{align*}
\begin{cases}
\beta u_{ttt}+u_{tt}-\Delta u-\beta\Delta u_t=\beta|u|^p,&x\in\mb{R}^n,\ t>0,\\
(u,u_t,u_{tt})(0,x)=(u_0,u_1,\Delta u_0)(x),&x\in\mb{R}^n.
\end{cases}
\end{align*}
Indeed, the deduction of the previous statement was motived by
\begin{align*}
\beta u_{ttt}+u_{tt}-\Delta u-\beta\Delta u_t&=(\beta\partial_t+\ml{I})(u_{tt}-\Delta u),\\
\beta|u|^p&=(\beta\partial_t+\ml{I})(g\ast|u|^p),
\end{align*}
where $\ml{I}$  is the identity operator such that $\ml{I}:f\to \ml{I}f\equiv f$. The recent paper \cite{ChenPalmieri201901} proved blow-up for the semilinear Moore-Gibson-Thompson equations in the conservative case if $1<p<p_{\mathrm{Str}}(n)$, which somehow coincides with our result in this example.
\end{exam}
\begin{exam}
Let us consider the semilinear wave equation \eqref{Eq Single semilinear wave Memory} with a very fast decay memory kernel such that
\begin{align*}
g(t)=\exp\left(\exp\left(\cdots\exp(-ct)\right)\right)\ \ \mbox{with}\ \ c>0.
\end{align*}
Due to the fact that $g(t)\lesssim t^{-1}$ for $t\geqslant t_0(c)$. We still can derive blow-up of solutions from Corollary \ref{Coro subcritical} if the power of the exponent fulfills $1<p<p_{\mathrm{Str}}(n)$.
\end{exam}
\subsection{Application of Theorem \ref{Thm subcritical}}
In this subsection, we will show several applications on Theorem \ref{Thm subcritical} by considering slow decay memory kernels and fast decay memory kernels, respectively. Throughout this part, we assume $p,q>1$ such that $p,q\leqslant n/(n-2)$ when $n\geqslant 3$.
\begin{exam}\label{Exam 01}
Let us consider the weakly coupled systems \eqref{Eq Semi Wave Distinct Memory} with nonlinear memory terms being the Riemann-Liouville fractional integrals of order $1-\gamma_1$ and $1-\gamma_2$, respectively, where $\gamma_1,\gamma_2\in(0,1)$. Precisely, we take
\begin{align}\label{Examp memory kernel}
g_1(t)=\frac{t^{-\gamma_1}}{\Gamma(1-\gamma_1)}\ \ \mbox{and}\ \  g_2(t)=\frac{t^{-\gamma_2}}{\Gamma(1-\gamma_2)}.
\end{align}
 Clearly, $g_k(t)\gtrsim t^{-1}$ and $g_k(t)\in L^1([0,T])$ for $k=1,2$ and any $t\geqslant 0$. We now employ Theorem \ref{Thm subcritical} with the condition \eqref{Hypothesis Condition 1} such that one of the following inequalities hold:
\begin{align*}
\frac{n-1}{2}pq-2q-\frac{1}{p}-\frac{n+5}{2}+\gamma_1q+\gamma_2<0\ \ \mbox{iff}\ \ \frac{(2-\gamma_1)q+(3-\gamma_2)+p^{-1}}{pq-1}>\frac{n-1}{2},\\
\frac{n-1}{2}pq-2p-\frac{1}{q}-\frac{n+5}{2}+\gamma_1+\gamma_2p<0\ \ \mbox{iff}\ \ \frac{(2-\gamma_2)p+(3-\gamma_1)+q^{-1}}{pq-1}>\frac{n-1}{2}.
\end{align*}
So, let us assume initial data satisfying the same assumptions as those in Theorem \ref{Thm subcritical} and $\alpha_{\mathrm{WM}}(p,q,\gamma_1,\gamma_2)>(n-1)/2$, where
\begin{align}\label{Examp Condition}
\alpha_{\mathrm{WM}}(p,q,\gamma_1,\gamma_2):=\max\left\{\frac{(2-\gamma_2)p+(3-\gamma_1)+q^{-1}}{pq-1},\frac{(2-\gamma_1)q+(3-\gamma_2)+p^{-1}}{pq-1}\right\}.
\end{align}
Then, every energy solution according to Definition \ref{Defn Energy Solution} of the weakly coupled system \eqref{Eq Semi Wave Distinct Memory} carrying \eqref{Examp memory kernel} blows up in finite time. Particularly, because the next relation holds in the sense of distribution:
\begin{align*}
\lim\limits_{\gamma\to 1^-}\frac{s_+^{-\gamma_k}}{\Gamma(1-\gamma_k)}=\delta_0(s)\ \ \mbox{carrying}\ \ s_+^{-\gamma_k}:=\begin{cases}
s^{-\gamma_k}&\mbox{if}\ \ s>0,\\
0&\mbox{if}\ \ s<0,
\end{cases}
\end{align*}
 it seems suitable to prove the blow-up result with \eqref{Examp Condition} satisfying
\begin{align*}
\lim\limits_{\gamma_1\to 1^-,\,\gamma_2\to 1^-}\alpha_{\mathrm{WM}}(p,q,\gamma_1,\gamma_2)=\alpha_{\mathrm{W}}(p,q),
\end{align*}
in which $\alpha_{\mathrm{W}}(p,q)$ describes the critical curve of the weakly coupled systems for wave equations \eqref{Eq Semi Wave Coupled}.
\end{exam}
\begin{exam}
Let us consider the weakly coupled systems \eqref{Eq Semi Wave Distinct Memory} with nonlinear memory terms, which are polynomial decay (faster than those in Example \ref{Exam 01}) without any singularities at $t=0$. To be specific, we consider
\begin{align*}
g_1(t)=(1+t)^{-\gamma_1}\ \ \mbox{and}\ \ g_2(t)=(1+t)^{-\gamma_2},
\end{align*}
where $\gamma_1,\gamma_2\in[1,\infty)$. By concerning the condition \eqref{Hypothesis Condition 2} and the same assumptions on initial data as those in Theorem \ref{Thm subcritical}, every energy solution $(u,v)$ blows up in finite time. Furthermore, if we consider the exponential decay memory kernels such that
\begin{align*}
g_1(t)=\mathrm{e}^{-c_1t}\ \ \mbox{and}\ \ g_2(t)=\mathrm{e}^{-c_2t},
\end{align*}
with $c_1,c_2>0$, then the blow-up condition is still shown by \eqref{Hypothesis Condition 2}. In the special case by following the same idea as Example \ref{Examp single 2}, we take $c_1=1/\beta_1>0$ and $c_2=1/\beta_2>0$, every energy solution to the weakly coupled systems for the Moore-Gibson-Thompson equations
\begin{align*}
\begin{cases}
\beta_1u_{ttt}+u_{tt}-\Delta u-\beta_1\Delta u_t=\beta_1|v|^p,&x\in\mb{R}^n,\ t>0,\\
\beta_2v_{ttt}+v_{tt}-\Delta v-\beta_2\Delta v_t=\beta_2|u|^q,&x\in\mb{R}^n,\ t>0,\\
(u,u_t,u_{tt},v,v_t,v_{tt})(0,x)=(u_0,u_1,\Delta u_0,v_0,v_1,\Delta v_0)(x),&x\in\mb{R}^n,
\end{cases}
\end{align*}
blows up provided that the condition \eqref{Hypothesis Condition 2} holds.
\end{exam}


\section{Proof of Theorem \ref{Thm Local existence}}\label{Section Proof of LOCAL EXISTENCE}
First of all, let us introduce some notations and well-developed results for linear wave equation. It is well-known that the solutions $\ml{w}=\ml{w}(t,x)$ of linear wave equation 
\begin{align}\label{Linear WAVE}
\begin{cases}
\ml{w}_{tt}-\Delta \ml{w}=0,&x\in\mb{R}^n,\,t>0,\\
(\ml{w},\ml{w}_t)(0,x)=(\ml{w}_0,\ml{w}_1)(x),&x\in\mb{R}^n,
\end{cases}
\end{align}
is given by the form
\begin{align*}
\ml{w}(t,x)=\ml{K}_0(t,x)\ast_{(x)}\ml{w}_0(x)+\ml{K}_1(t,x)\ast_{(x)}\ml{w}_1(x),
\end{align*}
where we denoted by $\ml{K}_0(t,x)$ and $\ml{K}_1(t,x)$ the fundamental solutions to the linear wave equation in $\mb{R}^n$ with initial data $(\ml{w}_0,\ml{w}_1)=(\delta_0,0)$ and $(\ml{w}_0,\ml{w}_1)=(0,\delta_0)$, respectively. Here, $\delta_0$ is the Dirac distribution in $x=0$ with respect to spatial variables. Precisely, the above kernels $\ml{K}_0(t,x)$ and $\ml{K}_1(t,x)$ may be represented by the application of partial Fourier transforms with respect to time variable as follows:
\begin{align*}
\ml{K}_0(t,x)=\ml{F}^{-1}_{\xi\to x}\left(\cos(|\xi|t)\right)\ \ \mbox{and}\ \ \ml{K}_1(t,x)=\ml{F}^{-1}_{\xi\to x}\left(\frac{\sin(|\xi|t)}{|\xi|}\right).
\end{align*}
 Moreover, the solution of the Cauchy problem \eqref{Linear WAVE} fulfills the following $L^2-L^2$ estimates:
\begin{align*}
\|\ml{w}(t,\cdot)\|_{L^2(\mb{R}^n)}&\lesssim\|\ml{w}_0\|_{L^2(\mb{R}^n)}+(1+t)\|\ml{w}_1\|_{L^2(\mb{R}^n)},\\
\|\nabla \ml{w}(t,\cdot)\|_{L^2(\mb{R}^n)}+\|\ml{w}_t(t,\cdot)\|_{L^2(\mb{R}^n)}&= \|\ml{w}_0\|_{H^1(\mb{R}^n)}+\|\ml{w}_1\|_{L^2(\mb{R}^n)}.
\end{align*}
These estimates can be proved easily by using the phase space analysis. One may check Chapter 10 of \cite{Ebert-Reissig-book} in detail.

We define the solution space to \eqref{Eq Semi Wave Distinct Memory} by
\begin{align*}
X(T):=&\big\{(u,v)\in\left(\ml{C}([0,T],H^1(\mb{R}^n))\cap \ml{C}^1([0,T],L^2(\mb{R}^n))\right)^2\ \ \mbox{such that}\\ &\quad\mathrm{supp}\,u(t,\cdot),\mathrm{supp}\,v(t,\cdot)\subset B_{R+t} \ \ \mbox{for any}\ \  t\in[0,T]\big\}
\end{align*}
carrying its norm
\begin{align*}
\|(u,v)\|_{X(T)}:=\max\limits_{t\in[0,T]}\left(\ml{M}[u](t)+\ml{M}[v](t)\right),
\end{align*}
where we defined for $w=u,v$ that
\begin{align*}
\ml{M}[w](t):=\|w(t,\cdot)\|_{L^2(\mb{R}^n)}+\|\nabla w(t,\cdot)\|_{L^2(\mb{R}^n)}+\|w_t(t,\cdot)\|_{L^2(\mb{R}^n)}.
\end{align*}
By using Duhamel's principle, we may introduce the operator $N$ such that
\begin{align*}
N:\ (u,v)\in X(T)\to N(u,v):=\left(u^{\lin}(t,x)+u^{\non}(t,x),v^{\lin}(t,x)+v^{\non}(t,x)\right).
\end{align*}
Here, we denote
\begin{align*}
u^{\lin}(t,x)&:=\ml{K}_0(t,x)\ast_{(x)}u_0(x)+\ml{K}_1(t,x)\ast_{(x)}u_1(x),\\
u^{\non}(t,x)&:=\int_0^t\ml{K}_1(t-\tau,x)\ast_{(x)}(g_1\ast|v|^p)(\tau,x)\mathrm{d}\tau,\\
v^{\lin}(t,x)&:=\ml{K}_0(t,x)\ast_{(x)}v_0(x)+\ml{K}_1(t,x)\ast_{(x)}v_1(x),\\
v^{\non}(t,x)&:=\int_0^t\ml{K}_1(t-\tau,x)\ast_{(x)}(g_2\ast|u|^q)(\tau,x)\mathrm{d}\tau.
\end{align*}
Next, we are going to consider as mild local (in time) solution of \eqref{Eq Semi Wave Distinct Memory} the fixed points of the operator $N$. In other words, we should prove
\begin{align}
\|N(u,v)\|_{X(T)}&\leqslant\bar{c}_0+\bar{c}_1T\left(\|(u,v)\|_{X(T)}^p+\|(u,v)\|_{X(T)}^q\right),\label{Cru 01}\\
\|N(u,v)-N(\bar{u},\bar{v})\|_{X(T)}&\leqslant\bar{c}_2T\|(u,v)-(\bar{u},\bar{v})\|_{X(T)}\sum\limits_{r=p,q}\left(\|(u,v)\|_{X(T)}^{r-1}+\|(\bar{u},\bar{v})\|_{X(T)}^{r-1}\right),\label{Cru 02}
\end{align}
where $\bar{c}_k=\bar{c}_k(u_0,u_1,v_0,v_1)>0$ for $k=1,2,3$, to showing the local (in time) existence and uniqueness of the solution in $X(T)$.

Obviously from the beginning of this section, we claim $(u^{\lin},v^{\lin})\in X(T)$ satisfying
\begin{align*}
\|(u^{\lin},v^{\lin})\|_{X(T)}&\lesssim\|(u_0,v_0)\|_{H^1(\mb{R}^n)\times H^1(\mb{R}^n)}+(1+T)\|(u_1,v_1)\|_{L^2(\mb{R}^n)\times L^2(\mb{R}^n)}.
\end{align*}

In order to arrive at \eqref{Cru 01}, by employing the classical Gagliardo-Nirenberg inequality, one immediately derives
\begin{align*}
\|w(\eta,\cdot)\|_{L^{2r}(\mb{R}^n)}^r\leqslant c\|w(\eta,\cdot)\|_{L^2(\mb{R}^n)}^{(1-\frac{n}{2}(1-\frac{1}{r}))r}\|\nabla w(\eta,\cdot)\|_{L^2(\mb{R}^n)}^{\frac{n}{2}(1-\frac{1}{r})r}\leqslant c\|(u,v)\|_{X(T)}^r,
\end{align*}
for $\eta\in[0,T]$, with $r>1$ if $n=1,2$ and $1<r\leqslant n/(n-2)$ if $n\geqslant 3$, where we denoted $w=u,v$ and $r=p,q$ in the last inequality.\\
Then, combining obtained $L^2-L^2$ estimates and the previous $L^{2r}$ estimates of solutions, we have
\begin{align*}
\left\|\int_0^t\ml{K}_1(t-\tau,x)\ast_{(x)}(g_1\ast|v|^p)(\tau,x)\right\|_{L^2(\mb{R}^n)}&\leqslant c\int_0^t(1+t-\tau)\int_0^{\tau}g_1(\tau-\eta)\|\,|v(\eta,\cdot)|^p\|_{L^2(\mb{R}^n)}\mathrm{d}\eta\mathrm{d}\tau\\
&\leqslant c\int_0^t(1+t-\tau)\int_0^{\tau}g_1(\tau-\eta)\mathrm{d}\eta\mathrm{d}\tau\,\|(u,v)\|_{X(T)}^p\\
&\leqslant c\int_0^t(1+t-\tau)\mathrm{d}\tau\,\|(u,v)\|_{X(T)}^p\\
&\leqslant c(1+t)t\|(u,v)\|_{X(T)}^p,
\end{align*}
where we used $g_1(t)\in L^1([0,T])$, namely, for $[0,\tau]\subset [0,T]$,
\begin{align*}
\int_0^{\tau}g_1(\tau-\eta)\mathrm{d}\eta=\int_0^{\tau}g_1(\eta)\mathrm{d}\eta\leqslant\int_0^Tg_1(t)\mathrm{d}t<\infty.
\end{align*}
Analogously,
\begin{align*}
\left\|\int_0^t\ml{K}_1(t-\tau,x)\ast_{(x)}(g_2\ast|u|^q)(\tau,x)\right\|_{L^2(\mb{R}^n)}\leqslant c(1+t)t\|(u,v)\|_{X(T)}^q,
\end{align*}
where we applied our assumption $g_2(t)\in L^1([0,T])$.\\
By the same way, we derived
\begin{align*}
\left\|\nabla^{j}\partial_t^k\int_0^t\ml{K}_1(t-\tau,x)\ast_{(x)}(g_1\ast|v|^p)(\tau,x)\mathrm{d}\tau\right\|_{L^2(\mb{R}^n)}\leqslant ct\|(u,v)\|^p_{X(T)},\\
\left\|\nabla^{j}\partial_t^k\int_0^t\ml{K}_1(t-\tau,x)\ast_{(x)}(g_2\ast|u|^q)(\tau,x)\mathrm{d}\tau\right\|_{L^2(\mb{R}^n)}\leqslant ct\|(u,v)\|^q_{X(T)},
\end{align*}
for any $j,k\in\mb{N}_0$ such that $j+k=1$. They implies our desired inequality \eqref{Cru 01}.\\
 Furthermore, the function $N(u,v)$ is the solution of the Cauchy problem for wave equations
\begin{align*}
\begin{cases}
\partial_t^2N(u,v)-\Delta N(u,v)=\left((g_1\ast|v|^p),(g_2\ast|u|^q)\right),&x\in\mb{R}^n,\,t>0,\\
(N(u,v),\partial_tN(u,v))(0,x)=(u_0,u_1,v_0,v_1)(x),&x\in\mb{R}^n.
\end{cases}
\end{align*}
We claim that $\mathrm{supp}\,N(u,v)(t,\cdot)\subset B_{R+t}\times B_{R+t}$ for any $t\in[0,T]$ because $(u,v)$ is assumed to be supported in the forward cone. Therefore, it is proved that the operator $N$ maps $X(T)$ into itself.

Finally, by considering $w=u,v$ and $r=p,q$, with the aid of
\begin{align*}
|\,|w(\tau,x)|^r-|\bar{w}(\tau,x)|^r|\leqslant c|w(\tau,x)-\bar{w}(\tau,x)|\left(|w(\tau,x)|^{r-1}+|\bar{w}(\tau,x)|^{r-1}\right),
\end{align*}
we are able to gain from H\"older's inequality that
\begin{align*}
\|\,|w(\tau,\cdot)|^r-|\bar{w}(\tau,\cdot)|^r\|_{L^2(\mb{R}^n)}\leqslant c\|w(\tau,\cdot)-\bar{w}(\tau,\cdot)\|_{L^2(\mb{R}^n)}\left(\|w(\tau,\cdot)\|_{L^{2r}(\mb{R}^n)}^{r-1}+\|\bar{w}(\tau,\cdot)\|_{L^{2r}(\mb{R}^n)}^{r-1}\right).
\end{align*}
We know the equality holds
\begin{align*}
\|N(u,v)-N(\bar{u},\bar{v})\|_{X(t)}=\left\|\int_0^t\ml{K}_1(t-\tau,x)\ast_{(x)}(F_p(\tau,x),F_q(\tau,x))\mathrm{d}\tau\right\|_{X(t)},
\end{align*}
where we denoted
\begin{align*}
F_p(\tau,x):=\int_0^{\tau}g_1(\tau-\eta)(|v(\eta,x)|^p-|\bar{v}(\eta,x)|^p)\mathrm{d}\eta,\\
F_q(\tau,x):=\int_0^{\tau}g_2(\tau-\eta)(|u(\eta,x)|^q-|\bar{u}(\eta,x)|^q)\mathrm{d}\eta.
\end{align*}
Similarly to the deduction of \eqref{Cru 01}, we use $L^2-L^2$ estimates, which allows us to conclude \eqref{Cru 02}. The proof is complete.

\section{Proof of Theorem \ref{Thm subcritical}}\label{Section Proof of THM 1}\setcounter{equation}{0}

To begin with the proof, let us first define time-dependent functionals related to the solutions
\begin{align}\label{Defn U(t) V(t)}
U(t):=\int_{\mb{R}^n}u(t,x)\mathrm{d}x\ \ \mbox{and}\ \ V(t):=\int_{\mb{R}^n}v(t,x)\mathrm{d}x.
\end{align}
In order to prove blow-up of energy solutions of the weakly coupled systems \eqref{Eq Semi Wave Distinct Memory}, we will show that the functional $(U(t),V(t))$ blows up in finite time by applying iteration methods. For this reason, we will construct a pair of coupled  of integral inequalities for the spatial averages of the components of a local (in time) solution $(u,v)$ by choosing special test functions in the definition of energy solutions, which is so-called iteration frame, and derive first lower bound estimates for two functionals. Then, sharper estimates for the functionals from the below will be established by employing suitable iteration argument.

Here, we take the test functions in \eqref{Eq Defn Energy 1} and \eqref{Eq Defn Energy 2} such that $\phi\equiv 1\equiv\psi$ in $\{(s,x)\in[0,t]\times\mb{R}^n:|x|\leqslant R+s\}$, which leads to
\begin{align*}
\int_{\mb{R}^n}u_t(t,x)\mathrm{d}x-\int_{\mb{R}^n}u_1(x)\mathrm{d}x&=\int_0^t\int_{\mb{R}^n}\int_0^sg_1(s-\tau)|v(\tau,x)|^p\mathrm{d}\tau\mathrm{d}x\mathrm{d}s,\\
\int_{\mb{R}^n}v_t(t,x)\mathrm{d}x-\int_{\mb{R}^n}v_1(x)\mathrm{d}x&=\int_0^t\int_{\mb{R}^n}\int_0^sg_2(s-\tau)|u(\tau,x)|^q\mathrm{d}\tau\mathrm{d}x\mathrm{d}s.
\end{align*}
The definition \eqref{Defn U(t) V(t)} allows us to rewrite the above relations into
\begin{align}
U'(t)&=U'(0)+\int_0^t\int_{\mb{R}^n}\int_0^sg_1(s-\tau)|v(\tau,x)|^p\mathrm{d}\tau\mathrm{d}x\mathrm{d}s,\label{Eq U'(t)}\\
V'(t)&=V'(0)+\int_0^t\int_{\mb{R}^n}\int_0^sg_2(s-\tau)|u(\tau,x)|^q\mathrm{d}\tau\mathrm{d}x\mathrm{d}s.\label{Eq V'(t)}
\end{align}
From the nonnegative assumptions on $u_0,u_1,v_0,v_1$ showing the nonnegativities of $U(0)$, $U'(0)$, $V(0)$, $V'(0)$, we are able to integrate \eqref{Eq U'(t)} and \eqref{Eq V'(t)}, respectively, over $[0,t]$ to obtain
\begin{align*}
U(t)&\geqslant\int_0^t\int_0^{\eta}\int_0^sg_1(s-\tau)\int_{\mb{R}^n}|v(\tau,x)|^p\mathrm{d}x\mathrm{d}\tau\mathrm{d}s\mathrm{d}\eta\geqslant0,\\
V(t)&\geqslant\int_0^t\int_0^{\eta}\int_0^sg_2(s-\tau)\int_{\mb{R}^n}|u(\tau,x)|^q\mathrm{d}x\mathrm{d}\tau\mathrm{d}s\mathrm{d}\eta\geqslant0.
\end{align*}
With the help of H\"older's inequality combined with the support conditions of solutions, there exist positive constants $C_0=C_0(n,R,p)$ and $\widetilde{C}_0=\widetilde{C}_0(n,R,q)$ such that
\begin{align*}
\int_{\mb{R}^n}|v(\tau,x)|^p\mathrm{d}x=\int_{B_{R+\tau}}|v(\tau,x)|^p\mathrm{d}x&\geqslant C_0(R+\tau)^{-n(p-1)}(V(\tau))^p,\\
\int_{\mb{R}^n}|u(\tau,x)|^q\mathrm{d}x=\int_{B_{R+\tau}}|u(\tau,x)|^q\mathrm{d}x&\geqslant \widetilde{C}_0(R+\tau)^{-n(q-1)}(U(\tau))^q.
\end{align*}
Therefore, the coupled system of crucial integral inequalities
\begin{align}
U(t)&\geqslant C_0\int_0^t\int_0^{\eta}\int_0^sg_1(s-\tau)(R+\tau)^{-n(p-1)}(V(\tau))^p\mathrm{d}\tau\mathrm{d}s\mathrm{d}\eta,\label{Iter Fra U}\\
V(t)&\geqslant \widetilde{C}_0\int_0^t\int_0^{\eta}\int_0^sg_2(s-\tau)(R+\tau)^{-n(q-1)}(U(\tau))^q\mathrm{d}\tau\mathrm{d}s\mathrm{d}\eta,\label{Iter Fra V}
\end{align}
hold for any $t\geqslant0$. In other words, the iteration frames  were constructed in \eqref{Iter Fra U} and \eqref{Iter Fra V}.

Our second step is to derive first lower bound estimates for the functionals $U(t)$ and $V(t)$, respectively. Strongly motivated by the pioneering research \cite{YordanovZhang2006}, we introduce the eigenfunction $\Phi=\Phi(x)$ of the Laplace operator in $n$ dimensional whole space, i.e.,
\begin{align*}
\Phi(x)&:=\mathrm{e}^{x}+\mathrm{e}^{-x}\qquad\qquad \ \mbox{if} \ \ n=1, \\  
\Phi(x)&:=\int_{\mathbb{S}^{n-1}}\mathrm{e}^{x\cdot \omega}\mathrm{d} \sigma_\omega \,\qquad \,\mbox{if} \ \  n\geqslant 2,
\end{align*}
where $\mathbb{S}^{n-1}$ is the $n-1$ dimensional sphere, which is devoted to the definition of the test function $\Psi=\Psi(t,x)$ such that $\Psi(t,x)=\mathrm{e}^{-t}\Phi(x)$. Plainly, the function $\Psi$ is the solution of the homogeneous wave equation $\Psi_{tt}-\Delta\Psi=0$. To investigate the estimates for $U(t)$ and $V(t)$ from the below, let us now define two auxiliary functionals
\begin{align*}
U_0(t):=\int_{\mb{R}^n}u(t,x)\Psi(t,x)\mathrm{d}x\ \ \mbox{and}\ \ V_0(t):=\int_{\mb{R}^n}v(t,x)\Psi(t,x)\mathrm{d}x.
\end{align*}
Differentiating \eqref{Eq U'(t)} and \eqref{Eq V'(t)} with respective to time variable and employing H\"older's inequality, compactness of the support of the solutions to have
\begin{align}
U''(t)&=\int_0^tg_1(t-\tau)\int_{\mb{R}^n}|v(\tau,x)|^p\mathrm{d}x\mathrm{d}\tau\geqslant C_1\int_0^tg_1(t-\tau)(R+\tau)^{n-1-\frac{n-1}{2}p}|V_0(\tau)|^p\mathrm{d}\tau\label{Est U''(t)},\\
V''(t)&=\int_0^tg_2(t-\tau)\int_{\mb{R}^n}|u(\tau,x)|^q\mathrm{d}x\mathrm{d}\tau\geqslant \widetilde{C}_1\int_0^tg_2(t-\tau)(R+\tau)^{n-1-\frac{n-1}{2}q}|U_0(\tau)|^q\mathrm{d}\tau\label{Est V''(t)},
\end{align}
with suitable positive constants $C_1=C_1(n,R,p)$ and $\widetilde{C}_1=\widetilde{C}_1(n,R,q)$. Here, we apply the asymptotic behavior of $\Psi$ to gain
\begin{align}\label{Est Behavior Psi(t,x)}
\int_{B_{R+\tau}}|\Psi(\tau,x)|^{\frac{r}{r-1}}\mathrm{d}x\lesssim (R+\tau)^{(n-1)(1-\frac{r'}{2})},
\end{align}
where $r'$ is the conjugate of $r$, i.e. $1/r+1/r'=1$. One also may find \eqref{Est Behavior Psi(t,x)} in \cite{LaiTakamura19}.\\
The nonnegativities of the nonlinearities of the equations in \eqref{Eq Semi Wave Distinct Memory} imply immediately from \cite{YordanovZhang2006} the lower bound of $U_0(t)$ and $V_0(t)$ that
\begin{align*}
U_0(t)&\geqslant\frac{1-\mathrm{e}^{-2t}}{2}\int_{\mb{R}^n}u_0(x)\Phi(x)\mathrm{d}x+\frac{1+\mathrm{e}^{-2t}}{2}\int_{\mb{R}^n}u_1(x)\Phi(x)\mathrm{d}x\geqslant C_2,\\
V_0(t)&\geqslant\frac{1-\mathrm{e}^{-2t}}{2}\int_{\mb{R}^n}v_0(x)\Phi(x)\mathrm{d}x+\frac{1+\mathrm{e}^{-2t}}{2}\int_{\mb{R}^n}v_1(x)\Phi(x)\mathrm{d}x\geqslant \widetilde{C}_2,
\end{align*}
for any $t\geqslant0$ with suitable constants $C_2>0$ and $\widetilde{C}_2>0$ depending on the size of initial data, where we used our assumptions on nontrivial data $u_1$ as well as $v_1$. In this step, one needs to employ the relations $\Psi_t(t,x)=-\Psi(t,x)$ and $\Delta\Psi(t,x)=\Psi(t,x)$.\\
Consequently, according to \eqref{Est U''(t)}, \eqref{Est V''(t)}, they lead to
\begin{align}
U(t)&\geqslant C_1\widetilde{C}_2^p\int_0^t\int_0^{\eta}\int_0^sg_1(s-\tau)(R+\tau)^{n-1-\frac{n-1}{2}p}\mathrm{d}\tau\mathrm{d}s\mathrm{d}\eta,\label{Sub Lower bound U(t)}\\
V(t)&\geqslant \widetilde{C}_1C_2^q\int_0^t\int_0^{\eta}\int_0^sg_2(s-\tau)(R+\tau)^{n-1-\frac{n-1}{2}q}\mathrm{d}\tau\mathrm{d}s\mathrm{d}\eta.\label{Sub Lower bound V(t)}
\end{align}
To understand the precise lower bound for the functionals, we need to discuss the estimates under different assumptions on the memory kernels. To be specific, we will complete the proof by considering the next two cases separately: 
\begin{itemize}
	\item Case 1: $g_k(t)\gtrsim t^{-1}$ for any $t\geqslant t_0$ with $t_0\in[0,T)$ and all $k=1,2$;
	\item Case 2: $g_k(t)\lesssim t^{-1}$ for any $t\geqslant t_0$ with $t_0\in[0,T)$ and all $k=1,2$.
\end{itemize}
This criterion $t^{-1}$ is motivated by the fact that if the time-dependent function $g_k(t)$ decreases very fast, then the following lower bound estimate:
\begin{align*}
\int_0^tg_k(t-\tau)\mathrm{d}\tau\gtrsim g_k(t)t\approx g_k(t)
\end{align*}
hold for any $t\gg 1$ and some losses appear, for example, the exponential decay memory kernel $g_k(t)=\mathrm{e}^{-t}$. It means that we will employ different approaches to derive lower bound estimates of the integral including the memory kernels, which primarily determined by the decreasing rate of the memory kernel.

\subsection{Treatment for Case 1}
In this case, we just need to directly apply the non-increasing properties for both memory kernels for any $t\geqslant0$ to catch from \eqref{Sub Lower bound U(t)} and \eqref{Sub Lower bound V(t)} that
\begin{align*}
U(t)&\geqslant\frac{C_1\widetilde{C}_2^p}{n(n+1)(n+2)}g_1(t)(R+t)^{-\frac{n-1}{2}p}t^{n+2},\\
V(t)&\geqslant\frac{\widetilde{C}_1C_2^q}{n(n+1)(n+2)}g_2(t)(R+t)^{-\frac{n-1}{2}q}t^{n+2}.
\end{align*}
\begin{remark}\label{Remark General}
Actually, one may generalize the assumption on $g_k(t)$ such that $g_k(t)\gtrsim \tilde{g}_k(t)>0$, where $\tilde{g}_k(t)$ is a constant or monotonously decreasing function. For example, from \eqref{Sub Lower bound U(t)} we may directly arrive at
\begin{align*}
U(t)&\gtrsim C_1\widetilde{C}_2^p(R+\tau)^{-\frac{n-1}{2}p}\int_0^t\int_0^{\eta}\int_0^s\tilde{g}_1(s-\tau)(R+\tau)^{n-1}\mathrm{d}\tau\mathrm{d}s\mathrm{d}\eta\\
&\gtrsim C_1\widetilde{C}_2^p\tilde{g}_1(t)(R+\tau)^{-\frac{n-1}{2}p}\int_0^t\int_0^{\eta}\int_0^s\tau^{n-1}\mathrm{d}\tau\mathrm{d}s\mathrm{d}\eta\\
&\gtrsim\frac{C_1\widetilde{C}_2^p}{n(n+1)(n+2)}\tilde{g}_1(t)(R+t)^{-\frac{n-1}{2}p}t^{n+2}.
\end{align*}
Then, we just need to replace $g_k(t)$ by $\tilde{g}_k(t)$ in the next all steps.
\end{remark}
We derived the first lower bounds for the functionals such that
\begin{align}
U(t)&\geqslant D_1(g_1(t))^{a_1}(g_2(t))^{\alpha_1}(R+t)^{-b_1}t^{\beta_1}\ \ \mbox{for any}\ \ t\geqslant 0,\label{First Lower Bound U}\\
V(t)&\geqslant \widetilde{D}_1(g_1(t))^{\tilde{a}_1}(g_2(t))^{\tilde{\alpha}_1}(R+t)^{-\tilde{b}_1}t^{\tilde{\beta}_1}\ \ \mbox{for any}\ \ t\geqslant 0.\label{First Lower Bound V}
\end{align}
Here, we denote constants in the previous estimates as
\begin{align*}
&D_1:=\frac{C_1\widetilde{C}_2^p}{n(n+1)(n+2)},\ \ \widetilde{D}_1:=\frac{\widetilde{C}_1C_2^q}{n(n+1)(n+2)},\\
&a_1:=1,\ \ \alpha_1:=0,\ \ b_1:=\frac{n-1}{2}p, \ \ \beta_1:=n+2,\\
&\tilde{a}_1:=0,\ \ \tilde{\alpha}_1:=1,\ \ \tilde{b}_1:=\frac{n-1}{2}q, \ \ \tilde{\beta}_1:=n+2,
\end{align*}

In the forthcoming part, by the way of iteration argument we will prove sequences of lower bounds of $U(t)$ and $V(t)$ by combining the iteration frame \eqref{Iter Fra U} and \eqref{Iter Fra V} as follows:
\begin{align}
U(t)&\geqslant D_j(g_1(t))^{a_j}(g_2(t))^{\alpha_j}(R+t)^{-b_j}t^{\beta_j}\ \ \mbox{for any}\ \ t\geqslant 0,\label{Sequence U}\\
V(t)&\geqslant \widetilde{D}_j(g_1(t))^{\tilde{a}_j}(g_2(t))^{\tilde{\alpha}_j}(R+t)^{-\tilde{b}_j}t^{\tilde{\beta}_j}\ \ \mbox{for any}\ \ t\geqslant 0,\label{Sequence V}
\end{align}
where $\{D_j\}_{j\geqslant1}$, $\{\widetilde{D}_j\}_{j\geqslant1}$, $\{a_j\}_{j\geqslant 1}$, $\{\tilde{a}_j\}_{j\geqslant 1}$, $\{\alpha_j\}_{j\geqslant1}$, $\{\tilde{\alpha}_j\}_{j\geqslant 1}$, $\{b_j\}_{j\geqslant 1}$, $\{\tilde{b}_j\}_{j\geqslant 1}$, $\{\beta_j\}_{j\geqslant 1}$ and $\{\tilde{\beta}_j\}_{j\geqslant 1}$ are sequences of nonnegative real numbers that will be determined iteratively in the inductive step. Note that the initial value, i.e. $j=0$, of the above inequalities are given in \eqref{First Lower Bound U} and \eqref{First Lower Bound V}. We should clarify that we cannot automatically ignore the components of $g_2(t)$ in the lower bound sequence for $U(t)$ and of $g_1(t)$ in the lower bound sequence for $V(t)$ although $\alpha_1=\tilde{a}_1=0$. The main factor comes from the interaction in weakly coupled systems.

For this reason, by assuming the validities of \eqref{Sequence U} and \eqref{Sequence V} for $j$, we just need to prove the induction step, i.e. the validities of \eqref{Sequence U} and \eqref{Sequence V} for $j+1$. Let us now combine \eqref{Sequence U}, \eqref{Sequence V} with \eqref{Iter Fra V}, \eqref{Iter Fra U} to arrive at
\begin{align*}
U(t)&\geqslant C_0\widetilde{D}_j^p\int_0^t\int_0^{\eta}\int_0^sg_1(s-\tau)(g_1(\tau))^{\tilde{a}_jp}(g_2(\tau))^{\tilde{\alpha}_jp}(R+\tau)^{-n(p-1)-\tilde{b}_jp}\tau^{\tilde{\beta}_jp}\mathrm{d}\tau\mathrm{d}s\mathrm{d}\eta\\
&\geqslant\frac{C_0\widetilde{D}_j^p}{(1+\tilde{\beta}_jp)(2+\tilde{\beta}_jp)(3+\tilde{\beta}_jp)}(g_1(t))^{1+\tilde{a}_jp}(g_2(t))^{\tilde{\alpha}_jp}(R+t)^{-n(p-1)-\tilde{b}_jp}t^{3+\tilde{\beta}_jp},
\end{align*}
and similarly,
\begin{align*}
V(t)&\geqslant\widetilde{C}_0D_j^q\int_0^t\int_0^{\eta}\int_0^s(g_1(\tau))^{a_jq}g_2(s-\tau)(g_2(\tau))^{\alpha_jq}(R+\tau)^{-n(q-1)-b_jq}\tau^{\beta_jq}\mathrm{d}\tau\mathrm{d}s\mathrm{d}\eta\\
&\geqslant\frac{\widetilde{C}_0D_j^q}{(1+\beta_jq)(2+\beta_jq)(3+\beta_jq)}(g_1(t))^{a_jq}(g_2(t))^{1+\alpha_jq}(R+t)^{-n(q-1)-b_jq}t^{3+\beta_jq}.
\end{align*}
Namely, the desired estimates \eqref{Sequence U} and \eqref{Sequence V} hold for $j+1$ providing that
\begin{align*}
&D_{j+1}:=\frac{C_0\widetilde{D}_j^p}{(1+\tilde{\beta}_jp)(2+\tilde{\beta}_jp)(3+\tilde{\beta}_jp)},\ \ \widetilde{D}_{j+1}:=\frac{\widetilde{C}_0D_j^q}{(1+\beta_jq)(2+\beta_jq)(3+\beta_jq)},\\
&a_{j+1}:=1+\tilde{a}_jp,\ \ \alpha_{j+1}:=\tilde{\alpha}_jp, \ \ b_{j+1}:=n(p-1)+\tilde{b}_jp, \ \  \beta_{j+1}:=3+\tilde{\beta}_jp,\\
& \tilde{a}_{j+1}:=a_jq,\ \  \tilde{\alpha}_{j+1}:=1+\alpha_jq,
 \ \ \tilde{b}_{j+1}:=n(q-1)+b_jq,\ \ \tilde{\beta}_{j+1}:=3+\beta_jq.
\end{align*}
In the proof, it is sufficient for us to determine the sequence for $a_j,\tilde{a}_j,\alpha_j,\tilde{\alpha}_j,b_j$ and $\tilde{b}_j$ for any odd number $j\geqslant3$. So, employing the last recursive relations with the initial value stated in \eqref{First Lower Bound U}, \eqref{First Lower Bound V} and the formula
\begin{align*}
a_j&=m+a_{j-2}pq=\cdots=m\left(1+(pq)+\cdots+(pq)^{\frac{j-3}{2}}\right)+a_1(pq)^{\frac{j-1}{2}}\\
&=\left(a_1+\frac{m}{pq-1}\right)(pq)^{\frac{j-1}{2}}-\frac{m}{pq-1},
\end{align*}
 one may get
\begin{align*}
&a_j=\frac{pq}{pq-1}(pq)^{\frac{j-1}{2}}-\frac{1}{pq-1},\ \ \,\,\,\, \tilde{a}_j=\frac{q}{pq-1}(pq)^{\frac{j-1}{2}}-\frac{q}{pq-1},\\
&\alpha_j=\frac{p}{pq-1}(pq)^{\frac{j-1}{2}}-\frac{p}{pq-1},\ \ \,\,\, \tilde{\alpha}_j=\frac{pq}{pq-1}(pq)^{\frac{j-1}{2}}-\frac{1}{pq-1},\\
&b_j=\frac{(n-1)p+2n}{2}(pq)^{\frac{j-1}{2}}-n,\ \ \tilde{b}_j=\frac{(n-1)q+2n}{2}(pq)^{\frac{j-1}{2}}-n,
\end{align*}
for any odd number $j$. By the same way, we deduce
\begin{align*}
\beta_j&=\frac{(n+2)(pq-1)+3(p+1)}{pq-1}(pq)^{\frac{j-1}{2}}-\frac{3(p+1)}{pq-1},\\
 \tilde{\beta}_j&=\frac{(n+2)(pq-1)+3(q+1)}{pq-1}(pq)^{\frac{j-1}{2}}-\frac{3(q+1)}{pq-1},
\end{align*}
for any odd number $j$, moreover, from the recursive relations between $\beta_j$ and $\tilde{\beta}_j$ again with even number $j$, i.e. $j-1$ is an odd number, one has 
\begin{align*}
\beta_j&=3+\tilde{\beta}_{j-1}p=\frac{(n+2)(pq-1)+3(q+1)}{(pq-1)q}(pq)^{\frac{j}{2}}-\frac{3(p+1)}{pq-1},\\
 \tilde{\beta}_j&=3+\beta_{j-1}q=\frac{(n+2)(pq-1)+3(p+1)}{(pq-1)p}(pq)^{\frac{j}{2}}-\frac{3(q+1)}{pq-1}.
\end{align*}
The previous representations of $\beta_j,\tilde{\beta}_j$ imply that there exist suitable positive constants $B_0$ and $\widetilde{B}_0$, which are independent of $j$, such that $\beta_j\leqslant B_0(pq)^{\frac{j}{2}}$ and $\tilde{\beta}_j\leqslant \widetilde{B}_0(pq)^{\frac{j}{2}}$ for any $j\geqslant1$. They immediately lead to
\begin{align*}
D_j&\geqslant C_0(2+\tilde{\beta}_{j-1}p)^{-3}\widetilde{D}_{j-1}^p\geqslant C_0\beta_j^{-3}\widetilde{D}_{j-1}^p\geqslant C_0B_0^{-3}(pq)^{-\frac{3}{2}j}\widetilde{D}_{j-1}^p,\\
\widetilde{D}_j&\geqslant \widetilde{C}_0(2+\beta_{j-1}q)^{-3}D_{j-1}^q\geqslant \widetilde{C}_0\tilde{\beta}_j^{-3}D_{j-1}^q\geqslant \widetilde{C}_0\widetilde{B}_0^{-3}(pq)^{-\frac{3}{2}j}D_{j-1}^q,
\end{align*}
for any odd number $j$. Summarizing the above relations yields
\begin{align*}
D_j&\geqslant C_0\widetilde{C}_0^pB_0^{-3}\widetilde{B}_0^{-3p}(pq)^{-\frac{3(p+1)}{2}j+\frac{3}{2}p}D_{j-2}^{pq}=:E_0(pq)^{-\frac{3(p+1)}{2}j}D_{j-2}^{pq},\\
\widetilde{D}_j&\geqslant C_0^q\widetilde{C}_0B_0^{-3q}\widetilde{B}_0^{-3}(pq)^{-\frac{3(q+1)}{2}j+\frac{3}{2}q}\widetilde{D}_{j-2}^{pq}=:\widetilde{E}_0(pq)^{-\frac{3(q+1)}{2}j}\widetilde{D}_{j-2}^{pq},
\end{align*}
for any odd number $j\geqslant 3$, where $E_0$ and $\widetilde{E}_0$ are suitable positive constants independent of $j$. Therefore, we may apply the logarithm on the recursive relationships to see 
\begin{align*}
\log D_j&\geqslant (pq)\log D_{j-2}-\frac{3(p+1)}{2}\log(pq)+\log E_0\\
&\geqslant\cdots\geqslant (pq)^{\frac{j-1}{2}}\log D_1-\frac{3(p+1)}{2}\log(pq)\sum\limits_{k=0}^{(j-3)/2}\left((j-2k)(pq)^k\right)+\log E_0\sum\limits_{k=0}^{(j-3)/2}(pq)^k\\
&=(pq)^{\frac{j-1}{2}}\left(\log D_1-\frac{3(p+1)(3pq-1)\log(pq)}{2(pq-1)^2}+\frac{\log E_0}{pq-1}\right)\\
&\quad+\frac{3(p+1)(2pq+j(pq-1))\log(pq)}{2(pq-1)^2}-\frac{\log E_0}{pq-1},
\end{align*}
and identically,
\begin{align*}
\log \widetilde{D}_j&\geqslant(pq)^{\frac{j-1}{2}}\left(\log \widetilde{D}_1-\frac{3(q+1)(3pq-1)\log(pq)}{2(pq-1)^2}+\frac{\log \widetilde{E}_0}{pq-1}\right)\\
&\quad+\frac{3(q+1)(2pq+j(pq-1))\log(pq)}{2(pq-1)^2}-\frac{\log \widetilde{E}_0}{pq-1},
\end{align*}
for any odd number $j\geqslant3$, where we applied the next formula:
\begin{align}\label{Sum formula}
\sum\limits_{k=0}^{(j-3)/2}\left((j-2k)(pq)^k\right)=\frac{2+3(pq-1)}{(pq-1)^2}(pq)^{\frac{j-1}{2}}-\frac{2pq+j(pq-1)}{(pq-1)^2}.
\end{align}
 Let us additionally consider the number $j$ fulfilling
\begin{align*}
j\geqslant j_0:=\left\lceil\frac{2}{3\log(pq)}\max\left\{\frac{\log E_0}{p+1},\frac{\log \widetilde{E}_0}{q+1}\right\}-\frac{2pq}{pq-1}\right\rceil.
\end{align*}
It means that for any odd number with $j\geqslant j_0$, we can estimate the multiplicative constants by
\begin{align*}
\log D_j&\geqslant(pq)^{\frac{j-1}{2}}\log\left(D_1(pq)^{-\frac{3(p+1)(3pq-1)}{2(pq-1)^2}}E_0^{\frac{1}{pq-1}}\right)=:(pq)^{\frac{j-1}{2}}\log E_1,\\
\log \widetilde{D}_j&\geqslant(pq)^{\frac{j-1}{2}}\log\left(\widetilde{D}_1(pq)^{-\frac{3(q+1)(3pq-1)}{2(pq-1)^2}}\widetilde{E}_0^{\frac{1}{pq-1}}\right)=:(pq)^{\frac{j-1}{2}}\log \widetilde{E}_1,
\end{align*}
where $E_1,\widetilde{E}_1$ are suitable positive constants independent of $j$.

In addition, concerning $t\geqslant R$ so that $\log(R+t)\leqslant\log(2t)$, we summarize the derived estimates and representations of constants in \eqref{Sequence U}, \eqref{Sequence V} to obtain
\begin{align}\label{Lower Bound U(t)}
U(t)&\geqslant\exp\left((pq)^{\frac{j-1}{2}}\log\left(E_12^{-\frac{(n-1)p}{2}-n}(g_1(t))^{\frac{pq}{pq-1}}(g_2(t))^{\frac{p}{pq-1}}t^{-\frac{(n-1)p}{2}+2+\frac{3(p+1)}{pq-1}}\right)\right)\notag\\
&\quad\times(g_1(t))^{-\frac{1}{pq-1}}(g_2(t))^{-\frac{p}{pq-1}}(R+t)^{n}t^{-\frac{3(p+1)}{pq-1}},
\end{align}
and
\begin{align}\label{Lower Bound V(t)}
V(t)&\geqslant\exp\left((pq)^{\frac{j-1}{2}}\log\left(\widetilde{E}_12^{-\frac{(n-1)q}{2}-n}(g_1(t))^{\frac{q}{pq-1}}(g_2(t))^{\frac{pq}{pq-1}}t^{-\frac{(n-1)q}{2}+2+\frac{3(q+1)}{pq-1}}\right)\right)\notag\\
&\quad\times (g_1(t))^{-\frac{q}{pq-1}}(g_2(t))^{-\frac{1}{pq-1}}(R+t)^{n}t^{-\frac{3(q+1)}{pq-1}},
\end{align}
for any odd number satisfying $j\geqslant j_0$. The assumption \eqref{Hypothesis Condition 1} is equivalent to one of the following conditions:
\begin{align*}
(g_1(t))^{\frac{pq}{pq-1}}(g_2(t))^{\frac{p}{pq-1}}t^{-\frac{(n-1)p}{2}+2+\frac{3(p+1)}{pq-1}}&\geqslant C(\ml{L}(t))^{\frac{p}{pq-1}},\\
(g_1(t))^{\frac{q}{pq-1}}(g_2(t))^{\frac{pq}{pq-1}}t^{-\frac{(n-1)q}{2}+2+\frac{3(q+1)}{pq-1}}&\geqslant C(\ml{L}(t))^{\frac{q}{pq-1}},
\end{align*}
hold. There exists a nonnegative constant $t_1$ such that for any $t\geqslant t_1$, it holds
\begin{align*}
(\ml{L}(t))^{\frac{\min\{p,q\}}{pq-1}}\geqslant(\ml{L}(t_1))^{\frac{\min\{p,q\}}{pq-1}}>\underbrace{C^{-1}\max\left\{E_1^{-1}2^{\frac{(n-1)p}{2}+n},\widetilde{E}_1^{-1}2^{\frac{(n-1)q}{2}+n}\right\}}_{\text{independent of}\,\,j}.
\end{align*}
Therefore, for $t\geqslant \max\{R,t_1\}$ and letting $j\to\infty$ in \eqref{Lower Bound U(t)} or \eqref{Lower Bound V(t)}, we may conclude that the lower bound for the functional $(U(t),V(t))$ blows up in finite time. Then, the proof of the theorem in Case 1 is complete.

\subsection{Treatment for Case 2}
Before discussing the blow-up result in Case 2, we set
\begin{align*}
G_k(t):=\int_0^tg_k(\tau)\mathrm{d}\tau\ \ \mbox{with}\ \ G'_k(t)=g_k(t)>0
\end{align*}
for $k=1,2$. Namely, $G_k(t)$ is an strictly increasing function for all $k=1,2$. Then, by the change of variable and employing integration by parts, one finds that
\begin{align}\label{Key Case 2}
\int_0^tg_k(t-\tau)\tau^{\alpha}\mathrm{d}\tau&=\int_0^tg_k(\tau)(t-\tau)^{\alpha}\mathrm{d}\tau=G_k(\tau)(t-\tau)^{\alpha}\big|_{\tau=0}^{\tau=t}+\alpha\int_0^tG_k(\tau)(t-\tau)^{\alpha-1}\mathrm{d}\tau\notag\\
&\geqslant
\begin{cases}
G_k(t_0)&\mbox{if}\ \ \alpha=0,\\
\displaystyle{\alpha\int_{t_0}^tG_k(\tau)(t-\tau)^{\alpha-1}\mathrm{d}}\tau&\mbox{if}\ \ \alpha>0,
\end{cases}\notag\\
&\geqslant  G_k(t_0)(t-t_0)^{\alpha}
\end{align}
for any $t\geqslant t_0$, where $\alpha\geqslant0$ and we used $G_k(0)=0$  and monotonous  properties of the function $G_k(t)$ for all $k=1,2$. Here, we shrank the domain of the integration from $[0,t]$ into $[t_0,t]$. Hence, the treatment \eqref{Key Case 2} motivate us to deal with the fast decay memory kernel in a better way.

We now can deduce the first lower bound estimates from \eqref{Sub Lower bound U(t)} that
\begin{align*}
U(t)&\geqslant C_1\widetilde{C}_2^p\int_0^t\int_0^\eta\int_0^sg_1(s-\tau)(R+\tau)^{-\frac{n-1}{2}p}\tau^{n-1}\mathrm{d}\tau\mathrm{d}s\mathrm{d}\eta\\
&\geqslant C_1\widetilde{C}_2^pG_1(t_0)(R+t)^{-\frac{n-1}{2}p}\int_{t_0}^t\int_{t_0}^{\eta}(s-t_0)^{n-1}\mathrm{d}s\mathrm{d}\eta\\
&\geqslant\frac{C_1\widetilde{C}_2^p}{n(n+1)}G_1(t_0)(R+t)^{-\frac{n-1}{2}p}(t-t_0)^{n+1}
\end{align*}
for any $t\geqslant t_0$, where in the first line of the chain inequality, the derived estimate \eqref{Key Case 2} was applied. Analogously, the estimate \eqref{Sub Lower bound V(t)} shows
\begin{align*}
V(t)\geqslant\frac{\widetilde{C}_1C_2^q}{n(n+1)}G_2(t_0)(R+t)^{-\frac{n-1}{2}q}(t-t_0)^{n+1}
\end{align*}
for any $t\geqslant t_0$. That is to say that we have derived the estimates for the functionals from the below as follows:
\begin{align}
U(t)\geqslant Q_1(R+t)^{-\theta_1}(t-L_1t_0)^{\sigma_1}\ \ \mbox{for any}\ \ t\geqslant L_1t_0,\label{Case2 Frist Lower U(t)}\\
V(t)\geqslant \widetilde{Q}_1(R+t)^{-\tilde{\theta}_1}(t-L_1t_0)^{\tilde{\sigma}_1}\ \ \mbox{for any}\ \ t\geqslant L_1t_0,\label{Case2 Frist Lower V(t)}
\end{align}
 where we set $L_1:=1$ and the constants are given by
\begin{align*}
&Q_1:=\frac{C_1\widetilde{C}_2^p}{n(n+1)}G_1(t_0),\ \ \widetilde{Q}_1:=\frac{\widetilde{C}_1C_2^q}{n(n+1)}G_2(t_0),\\
&\theta_1:=\frac{n-1}{2}p,\ \ \sigma_1:=n+1, \ \ \tilde{\theta}_1:=\frac{n-1}{2}q,\ \ \tilde{\sigma}_1:=n+1.
\end{align*}

Similarly to the previous subsection, we will demonstrate the sequence of lower bounds for $U(t)$ and $V(t)$ by the next way:
\begin{align}
U(t)\geqslant Q_j(R+t)^{-\theta_j}(t-L_jt_0)^{\sigma_j}\ \ \mbox{for any}\ \ t\geqslant L_jt_0,\label{Case2 Seq U(t)}\\
V(t)\geqslant \widetilde{Q}_j(R+t)^{-\tilde{\theta}_j}(t-L_jt_0)^{\tilde{\sigma}_j}\ \ \mbox{for any}\ \ t\geqslant L_jt_0,\label{Case2 Seq V(t)}
\end{align}
 where $\{Q_j\}_{j\geqslant 1}$, $\{\widetilde{Q}_j\}_{j\geqslant 1}$, $\{\theta_j\}_{j\geqslant 1}$, $\{\tilde{\theta}_j\}_{j\geqslant 1}$, $\{\sigma\}_{j\geqslant 1}$ and $\{\tilde{\sigma}\}_{j\geqslant 1}$ are sequences of nonnegative real numbers that will be determined later. Furthermore, motivated by the recent paper \cite{ChenPalmieri201901}, we construct $\{L_j\}_{j\geqslant 1}$ to be the sequence of the partial products of the convergent infinite product
\begin{align}\label{Seq_Nakao}
\prod\limits_{k=1}^{\infty}\ell_{k}\ \ \mbox{with}\ \ \ell_k:=1+(pq)^{-\frac{k-1}{2}}\ \ \mbox{for any}\ \  k\geqslant1,
\end{align}
that is,
\begin{align}\label{Seq_Nakao_2}
L_j:=\prod\limits_{k=1}^{j}\ell_{k}\ \ \mbox{for any}\ \  j\geqslant1.
\end{align}
Clearly, the convergent property can be easily obtained from the ratio test method and the fact that $\lim_{k\to\infty}(\ln\ell_{k+1})/(\ln\ell_k)=(pq)^{-1/2}<1$. Note that the initial value of the previous sequences was stated in \eqref{Case2 Frist Lower U(t)} as well as \eqref{Case2 Frist Lower V(t)}. 

Let us begin with processing the inductive step by combining \eqref{Case2 Seq U(t)}, \eqref{Case2 Seq V(t)} and \eqref{Iter Fra U}, \eqref{Iter Fra V}. We assume \eqref{Case2 Seq U(t)} and \eqref{Case2 Seq V(t)} hold for $j$, and our goal is to prove them also hold for $j+1$. By this way, one may get
\begin{align*}
U(t)\geqslant C_0\widetilde{Q}_j^p(R+t)^{-n(p-1)-\tilde{\theta}_jp}\int_{L_jt_0}^t\int_{L_jt_0}^{\eta}\int_{L_jt_0}^sg_1(s-\tau)(\tau-L_jt_0)^{\tilde{\sigma}_jp}\mathrm{d}\tau\mathrm{d}s\mathrm{d}\eta.
\end{align*}
Let us discuss the estimate the $\tau$-integral in the previous line. By changing of the variable and using integration by parts associated with $G_1(0)=0$, we have
\begin{align*}
\int_{L_jt_0}^sg_1(s-\tau)(\tau-L_jt_0)^{\tilde{\sigma}_jp}\mathrm{d}\tau&=\int_0^{s-L_jt_0}g_1(\tau)(s-L_jt_0-\tau)^{\tilde{\sigma}_jp}\mathrm{d}\tau\\
&=\tilde{\sigma}_jp\int_0^{s-L_jt_0}G_1(\tau)(s-L_jt_0-\tau)^{\tilde{\sigma}_jp-1}\mathrm{d}\tau\\
&\geqslant\tilde{\sigma}_jp\int_{L_jt_0(\ell_{j+1}-1)}^{s-L_jt_0}G_1(\tau)(s-L_jt_0-\tau)^{\tilde{\sigma}_jp-1}\mathrm{d}\tau\\
&\geqslant G_1(L_jt_0(\ell_{j+1}-1))(s-L_jt_0\ell_{j+1})^{\tilde{\sigma}_jp}
\end{align*}
for any $s\geqslant L_{j+1}t_0$, where we shrank the domain by using the fact that
\begin{align*}
s\geqslant L_{j+1}t_0=L_jt_0\ell_{j+1}\ \ \mbox{implies}\ \ s-L_jt_0\geqslant L_jt_0(\ell_{j+1}-1)>0.
\end{align*}
Let us estimate $G_1(L_jt_0(\ell_{j+1}-1))$ from below now by introducing
\begin{align}\label{Limt}
L:= \lim\limits_{j\to\infty}L_j=\prod\limits_{j=1}^{\infty}\ell_j=\prod\limits_{j=1}^{\infty}\left(1+(pq)^{-\frac{j-1}{2}}\right)>1.
\end{align}
Since $\ell_j>1$ the sequence $\{L_j\}_{j\geqslant1}$ is converging to $L$ as $j\to\infty$. So, we may claim $1\leqslant L_j\leqslant L$ for all $j\geqslant 1$. There is a positive real integer $j_m$ such that if $j\geqslant j_m$, then it holds
\begin{align*}
0<L_jt_0(\ell_{j+1}-1)\leqslant Lt_0(pq)^{-j/2}\ll 1.
\end{align*}
At this time, according to our assumptions $g_1(t)\in\ml{C}^2([0,T])$ and $g''_1(0)>0$, it would lead to
\begin{align*}
G_1(L_jt_0(\ell_{j+1}-1))&\geqslant G_1(0)+g_1(0)L_jt_0(pq)^{-\frac{j}{2}}+\frac{g'_1(0)}{2}L_j^2t_0^2(pq)^{-j}+g''_1(0)\ml{O}\left((pq)^{-\frac{3}{2}j}L_j^3\right)\\
&\geqslant(pq)^{-j}L_jt_0\left((pq)^{\frac{j}{2}}g_1(0)+\frac{g'_1(0)}{2}L_jt_0\right)+g''_1(0)\ml{O}\left((pq)^{-\frac{3}{2}j}\right)\\
&\geqslant(pq)^{-j}t_0\left((pq)^{\frac{j}{2}}g_1(0)-\frac{-g'_1(0)}{2}Lt_0\right)\\
&\geqslant C_3(pq)^{-j}
\end{align*}
for any $j\geqslant \max\{j_m,j_1\}$, where $C_3>0$ and we denoted
\begin{align*}
j_1:=\begin{cases}
1&\mbox{if}\ \ g_1'(0)>0,\\
\left\lceil2\log_{pq}\left(\frac{1}{g_1(0)}-\frac{g_1'(0)Lt_0}{2g_1(0)}\right)\right\rceil&\mbox{if}\ \ g_1'(0)\leqslant0.
\end{cases}
\end{align*}
We should emphasize that the definition of $j_1$ when $g'_1(0)\leqslant 0$ turn out from
\begin{align*}
(pq)^{\frac{j}{2}}g_1(0)-\frac{-g'_1(0)}{2}Lt_0\geqslant 1
\end{align*}
for any integer $j\geqslant j_0$\\
Summing up the derived estimates, we conclude
\begin{align*}
U(t)&\geqslant C_0C_3(pq)^{-j}\widetilde{Q}_j^p(R+t)^{-n(p-1)-\tilde{\theta}_jp}\int_{L_{j+1}t_0}^t\int_{L_{j+1}t_0}^{\eta}(s-L_jt_0\ell_{j+1})^{\tilde{\sigma}_jp}\mathrm{d}s\mathrm{d}\eta\\
&\geqslant\frac{C_0C_3(pq)^{-j}\widetilde{Q}_j^p}{(\tilde{\sigma}_jp+1)(\tilde{\sigma}_jp+2)}(R+t)^{-n(p-1)-\tilde{\theta}_jp}(t-L_{j+1}t_0)^{\tilde{\sigma}_jp+2}
\end{align*}
for any $t\geqslant L_{j+1}t_0$ and $j\geqslant \max\{j_m,j_1\}$.\\
By the same way of deduction as the lower bound estimates for $U(t)$, there exists a positive constant $\widetilde{C}_3$ such that
\begin{align*}
V(t)\geqslant\frac{\widetilde{C}_0\widetilde{C}_3(pq)^{-j}Q_j^q}{(\sigma_jq+1)(\sigma_jq+2)}(R+t)^{-n(q-1)-\theta_jq}(t-L_{j+1}t_0)^{\sigma_jq+2}
\end{align*}
for any $t\geqslant L_{j+1}t_0$ and $j\geqslant \max\{j_m,\tilde{j}_1\}$ with
\begin{align*}
\tilde{j}_1:=\begin{cases}
1&\mbox{if}\ \ g_2'(0)>0,\\
\left\lceil2\log_{pq}\left(\frac{1}{g_2(0)}-\frac{g_2'(0)Lt_0}{2g_2(0)}\right)\right\rceil&\mbox{if}\ \ g_2'(0)\leqslant0.
\end{cases}
\end{align*}
In other words, the desired lower bound estimates \eqref{Case2 Seq U(t)} and \eqref{Case2 Seq V(t)} are valid for $j+1$ if
\begin{align*}
&Q_{j+1}:=\frac{C_0C_3(pq)^{-j}\widetilde{Q}_j^p}{(\tilde{\sigma}_jp+1)(\tilde{\sigma}_jp+2)},\ \ \widetilde{Q}_{j+1}:=\frac{\widetilde{C}_0\widetilde{C}_3(pq)^{-j}Q_j^q}{(\sigma_jq+1)(\sigma_jq+2)},\\
&\theta_{j+1}:=n(p-1)+\tilde{\theta}_jp,\ \ \sigma_{j+1}:=\tilde{\sigma}_jp+2,\ \ \tilde{\theta}_{j+1}:=n(q-1)+\theta_jq,\ \ \tilde{\sigma}_{j+1}:=\sigma_jq+2.
\end{align*}

Similarly to the treatment in the last subsection, for any odd integer $j$ such that $j\geqslant\max\{j_m,j_1,\tilde{j}_1\}$, we employ iteratively the obtained relations to get
\begin{align*}
\theta_{j}=\frac{2n+(n-1)p}{2}(pq)^{\frac{j-1}{2}}-n, \ \ \tilde{\theta}_{j}=\frac{2n+(n-1)q}{2}(pq)^{\frac{j-1}{2}}-n.
\end{align*}
What's more, we obtain
\begin{align*}
\sigma_j&=\frac{(n+1)(pq-1)+2(p+1)}{pq-1}(pq)^{\frac{j-1}{2}}-\frac{2(p+1)}{pq-1},\\
\tilde{\sigma}_j&=\frac{(n+1)(pq-1)+2(q+1)}{pq-1}(pq)^{\frac{j-1}{2}}-\frac{2(q+1)}{pq-1},
\end{align*}
for any odd number $j\geqslant\max\{j_m,j_1,\tilde{j}_1\}$, and 
\begin{align*}
\sigma_j&=\tilde{\sigma}_{j-1}p+2=\frac{(n+1)(pq-1)+2(q+1)}{(pq-1)q}(pq)^{\frac{j}{2}}-\frac{2(p+1)}{pq-1},\\
\tilde{\sigma}_j&=\sigma_{j-1}q+2=\frac{(n+1)(pq-1)+2(q+1)}{(pq-1)p}(pq)^{\frac{j}{2}}-\frac{2(q+1)}{pq-1},
\end{align*}
for any even number $j\geqslant\max\{j_m,j_1,\tilde{j}_1\}$. It immediately leads to $\sigma_j\leqslant B_1(pq)^{\frac{j}{2}}$ and $\tilde{\sigma}_j\leqslant \widetilde{B}_1(pq)^{\frac{j}{2}}$ with suitable positive constants $B_1$ and $\widetilde{B}_1$ independent of $j$. 

Next, by using the relation between $Q_j$ and $\widetilde{Q}_j$ such that
\begin{align*}
Q_j&\geqslant\frac{C_0C_3}{B_1^2}(pq)^{-2j+1}\widetilde{Q}_{j-1}^p\geqslant\frac{C_0\widetilde{C}_0^pC_3\widetilde{C}_3^{p}(pq)^{3p+1}}{B_1^2\widetilde{B}_1^{2p}}(pq)^{-2j(p+1)}Q_{j-2}^{pq}=:E_2(pq)^{-2j(p+1)}Q_{j-2}^{pq},\\
 \widetilde{Q}_j&\geqslant\frac{\widetilde{C}_0\widetilde{C}_3}{\widetilde{B}_1^2}(pq)^{-2j+1}Q_{j-1}^q\geqslant\frac{\widetilde{C}_0C_0^qC_3^q\widetilde{C}_3(pq)^{3q+1}}{B_1^{2q}\widetilde{B}_1^2}(pq)^{-2j(q+1)}\widetilde{Q}_{j-2}^{pq}=:\widetilde{E}_2(pq)^{-2j(q+1)}\widetilde{Q}_{j-2}^{pq},
\end{align*}
where $E_2$ and $\widetilde{E}_2$ are suitable positive constants, one may obtain by applying the logarithmic function on the both sides of them to have
\begin{align*}
\log Q_j&\geqslant (pq)^{\frac{j-1}{2}}\log Q_1-2(p+1)\log(pq)\sum\limits_{k=0}^{(j-3)/2}\left((j-2k)(pq)^k\right)+\log E_2\sum\limits_{k=0}^{(j-3)/2}(pq)^k\\
&=(pq)^{\frac{j-1}{2}}\left(\log Q_1-\frac{2(p+1)(3pq-1)\log(pq)}{(pq-1)^2}+\frac{\log E_2}{pq-1}\right)\\
&\quad+\frac{2(p+1)(2pq+j(pq-1))\log(pq)}{(pq-1)^2}-\frac{\log E_2}{pq-1},
\end{align*}
as well as
\begin{align*}
\log \widetilde{Q}_j&\geqslant (pq)^{\frac{j-1}{2}}\log \widetilde{Q}_1-\frac{2(q+1)}{2}\log(pq)\sum\limits_{k=0}^{(j-3)/2}\left((j-2k)(pq)^k\right)+\log \widetilde{E}_2\sum\limits_{k=0}^{(j-3)/2}(pq)^k\\
&=(pq)^{\frac{j-1}{2}}\left(\log \widetilde{Q}_1-\frac{2(q+1)(3pq-1)\log(pq)}{(pq-1)^2}+\frac{\log \widetilde{E}_2}{pq-1}\right)\\
&\quad+\frac{2(q+1)(2pq+j(pq-1))\log(pq)}{(pq-1)^2}-\frac{\log \widetilde{E}_2}{pq-1},
\end{align*}
for any odd number $j\geqslant\max\{j_m,j_1,\tilde{j}_1\}$, where we used \eqref{Sum formula} again. Fixing
\begin{align*}
j_2:=\left\lceil\frac{1}{2\log(pq)}\max\left\{\frac{\log E_2}{p+1},\frac{\log \widetilde{E}_2}{q+1}\right\}-\frac{2pq}{pq-1}\right\rceil,
\end{align*}
therefore, for any odd number $j\geqslant\max\{j_m,j_1,\tilde{j_1},j_2\}$, we can estimate the multiplicative constants as follows:
\begin{align*}
\log Q_j&\geqslant(pq)^{\frac{j-1}{2}}\log\left(Q_1(pq)^{-\frac{2(p+1)(3pq-1)}{(pq-1)^2}}E_2^{\frac{1}{pq-1}}\right)=:(pq)^{\frac{j-1}{2}}\log E_3,\\
\log \widetilde{Q}_j&\geqslant(pq)^{\frac{j-1}{2}}\log\left(\widetilde{Q}_1(pq)^{-\frac{2(q+1)(3pq-1)}{(pq-1)^2}}\widetilde{E}_2^{\frac{1}{pq-1}}\right)=:(pq)^{\frac{j-1}{2}}\log \widetilde{E}_3,
\end{align*}
where $E_3,\widetilde{E}_3$ are suitable positive constants independent of $j$.

Let us recall the limit \eqref{Limt}, which leads to \eqref{Case2 Seq U(t)} and \eqref{Case2 Seq V(t)}  for $t\geqslant Lt_0$. Concerning $t\geqslant \max\{R,2Lt_0\}$ showing
\begin{align*}
\log(R+t)\leqslant\log(2t)\ \ \mbox{and}\ \ \log(t-Lt_0)\geqslant\log(t/2),
\end{align*}
 we summarize the derived estimates and representations of constants in \eqref{Case2 Seq U(t)}, \eqref{Case2 Seq V(t)} to obtain
\begin{align}
U(t)&\geqslant\exp\left((pq)^{\frac{j-1}{2}}\log\left(Q_12^{-\frac{2n+(n-1)p}{2}-\frac{(n+1)(pq-1)+2(p+1)}{pq-1}}t^{-\frac{(n-1)p}{2}+1+\frac{2(p+1)}{pq-1}}\right)\right)(R+t)^{n}t^{-\frac{2(p+1)}{pq-1}},\label{Case 2 Lower Bound U(t)}\\
V(t)&\geqslant\exp\left((pq)^{\frac{j-1}{2}}\log\left(\widetilde{Q}_12^{-\frac{2n+(n-1)q}{2}-\frac{(n+1)(pq-1)+2(q+1)}{pq-1}}t^{-\frac{(n-1)q}{2}+1+\frac{2(q+1)}{pq-1}}\right)\right)(R+t)^{n}t^{-\frac{2(q+1)}{pq-1}},\label{Case Lower Bound V(t)}
\end{align}
for any odd number fulfilling $j\geqslant \max\{j_m,j_1,\tilde{j_1},j_2\}$.\\
 Let us recall the condition \eqref{Hypothesis Condition 2} which is equivalent to
\begin{align*}
-\frac{(n-1)p}{2}+1+\frac{2(p+1)}{pq-1}>0\ \ \mbox{or}\ \ -\frac{(n-1)q}{2}+1+\frac{2(q+1)}{pq-1}>0.
\end{align*}
The power of $t$ in the exponential function is positive. Eventually, by taking  $t\geqslant\{R,2Lt_0\}$ and
\begin{align*}
t>\min\left\{\left(Q_12^{2n+1+\frac{n-1}{2}p+\frac{2(p+1)}{pq-1}}\right)^{-\frac{(n-1)p}{2}+1+\frac{2(p+1)}{pq-1}},\left(\widetilde{Q}_12^{2n+1+\frac{n-1}{2}q+\frac{2(q+1)}{pq-1}}\right)^{-\frac{(n-1)q}{2}+1+\frac{2(q+1)}{pq-1}}\right\},
\end{align*}
we may immediately find blow-up for the lower bounds of the functional $(U(t),V(t))$ in \eqref{Case 2 Lower Bound U(t)} and \eqref{Case Lower Bound V(t)} as $j\to\infty$. In conclusion, this concludes the proof in Case 2.

\section{Final remarks}\label{Section final remark}\setcounter{equation}{0}

Throughout this paper, we have mainly derived blow-up of energy solutions for the weakly coupled systems \eqref{Eq Semi Wave Distinct Memory} for distinct memory terms in general forms. Precisely, we determine a threshold $t^{-1}$ for the blow-up condition. However, it is still open to determine the critical condition or the critical curve for \eqref{Eq Semi Wave Distinct Memory} under certain assumptions on memory kernels. 

Let us now show some open problems strongly related to the topic in this paper and give some conjectures on them. In the following discussion, we generally assume $g_k(t)\in L^1([0,T])$ with $T\leqslant \infty$ for $k=1,2$.
\begin{itemize}
	\item If $g_k(t)\gtrsim t^{-1}$ for any $t\geqslant t_0$ with $t_0\in[0,T)$ and $g'_k(t)\leqslant0$ for $k=1,2$, then it is still open for the existence of global (in time) solution in the supercritical case
	\begin{align*}
	g_1(t)g_2(t)\max\left\{(g_1(t))^{q-1}t^{2q+1/p},(g_2(t))^{p-1}t^{2p+1/q}\right\}\lesssim t^{\frac{n-1}{2}(pq-1)-3},
	\end{align*}
	and nonexistence of local (in time) solution in the critical case
	\begin{align*}
	g_1(t)g_2(t)\max\left\{(g_1(t))^{q-1}t^{2q+1/p},(g_2(t))^{p-1}t^{2p+1/q}\right\}\asymp t^{\frac{n-1}{2}(pq-1)-3}.
	\end{align*}
	We should underline that the conjecture in the critical case for $p=q$ and $g_1(t)=g_2(t)=t^{-\gamma}$ with $\gamma\in(0,1)$, from point of view of the single semilinear equation, is true and we refer to Theorem 2 in \cite{ChenPalmieri20}. Hence, the conjectures seem reasonable.
	\item If $g_k(t)\lesssim t^{-1}$ for any $t\geqslant t_0$ with $t_0\in[0,T)$ and $g_k(t)\in\ml{C}^2([0,T])$ with $g''_k(0)>0$ for $k=1,2$, then it is still open for the existence of global (in time) solution in the supercritical case
	\begin{align*}
	\max\left\{\frac{p+2+q^{-1}}{pq-1},\frac{q+2+p^{-1}}{pq-1}\right\}<\frac{n-1}{2},
	\end{align*}
	and nonexistence of local (in time) solution in the critical case
	\begin{align*}
	\max\left\{\frac{p+2+q^{-1}}{pq-1},\frac{q+2+p^{-1}}{pq-1}\right\}=\frac{n-1}{2}.
	\end{align*}
	We expect the nonexistence in the critical case can be proved by constructing a weighted space average as the functional, whose dynamic can be studied in the iteration procedure \cite{Agemi-Kurokawa-Takamura2000,ChenPalmieri201901}. By this way, the conjectures seem reasonable too.
\end{itemize}

\section*{Acknowledgments}
The Ph.D. study of Wenhui Chen are supported by S\"achsisches Landesgraduiertenstipendium. This work was partially written while Wenhui Chen was a Ph.D. student at TU Freiberg. The author thanks  Michael Reissig (TU Bergakademie Freiberg), Alessandro Palmieri (University of Pisa) and Giovanni Girardi (University of Bari) for the suggestions in the preparation of the paper.

\end{document}